\documentclass[a4paper,twoside]{article}
\usepackage[margin=1.8cm,bottom=1cm,includefoot]{geometry}
\usepackage{amsmath,amssymb,amsthm}
\usepackage{bbm}
\usepackage{array}
\usepackage{mathrsfs}
\usepackage[utf8]{inputenc}
\usepackage{hyperref}
\usepackage[all]{hypcap}
\usepackage{float}

\newtheorem{thm}{Theorem}
\newtheorem{lem}[thm]{Lemma}
\newtheorem{prop}[thm]{Proposition}

\let\BFseries\bfseries\def\bfseries{\BFseries\mathversion{bold}}

\newcommand{\N}{\mathbb{N}}
\newcommand{\Z}{\mathbb{Z}}
\newcommand{\Q}{\mathbb{Q}}
\newcommand{\R}{\mathbb{R}}

\newcommand{\E}{\mathbb{E}}
\newcommand{\p}{\mathbb{P}}
\newcommand{\CC}{\mathcal{C}}
\newcommand{\LL}{\mathcal{L}}
\newcommand{\BB}{\mathcal{B}}
\newcommand{\MM}{\mathcal{M}}
\newcommand{\UU}{\mathcal{U}}
\newcommand{\DD}{\mathcal{D}}
\newcommand{\EE}{\mathcal{E}}
\newcommand{\1}{\mathbbm{1}}
\newcommand{\eps}{\varepsilon}
\DeclareMathOperator{\e}{e}
\newcommand{\D}{\mathrm{d}}
\newcommand{\op}{\operatorname}
\newcommand{\lne}{<}
\newcommand{\gne}{>}
\newcommand{\abs}{|}
\newcommand{\y}{\|}
\newcommand{\calC}{{\mathcal C}}
\newcommand{\calD}{{\mathcal D}}
\newcommand{\pl}{\partial}

\newcommand{\iy}{\infty}

\newcommand{\Title}{Trimmed Branching Brownian Motion with Discontinuous Drift}
\newcommand{\Author}{Rami Atar, Leonid Mytnik, Dominic T. Schickentanz}

\makeatletter
\renewcommand{\@oddhead}{\scshape \hfill \Title \hfill}
\renewcommand{\@oddfoot}{\hfill \thepage \hfill}
\renewcommand{\@evenhead}{\scshape \hfill \Author \hfill}
\renewcommand{\@evenfoot}{\hfill \thepage \hfill}
\makeatother

\begin{document}
	\title{\Title}
	\author{Rami Atar\thanks{Viterbi Faculty of Electrical \& Computer Engineering, Technion. Email address: \texttt{rami@technion.ac.il}},\quad Leonid Mytnik\thanks{Faculty of Data \& Decision Sciences, Technion. Email address: \texttt{leonidm@technion.ac.il}},\quad Dominic T. Schickentanz\thanks{Faculty of Mathematics, Technion. Email address: \texttt{mail@dominic-schickentanz.com}}}
	\date{November 19, 2026}
	\maketitle

	\begin{abstract}
		We consider a system of branching Brownian particles on the real line with bounded, possibly discontinuous drift $b$, and with a selection rule that removes one particle from the most crowded cells of a uniform partition of the line at each branching event so that the number of particles is preserved.
		As the number of particles tends to infinity, with the cell size tending to zero (sufficiently slowly), we prove convergence in probability of the empirical measure process to a deterministic limit characterized by the parabolic free obstacle problem $$\partial_tu=\tfrac12\partial_{xx}u-\partial_x(bu)+u-\beta.$$
		Here the density $u$ has unit mass. Disintegrating the removal measure $\beta$ as $\beta_t(\mathrm{d} x)\mathrm{d} t$, the probability measure~$\beta_t$ is supported on the spatial maximizers of $u(t,\cdot)$ for a.e.~$t$. We establish existence and uniqueness of weak solutions in a suitable class of continuous densities, extending earlier results about this PDE to discontinuous drifts.
		We then consider an inverse problem: Given a target removal measure, find a drift for which the PDE admits a stationary solution with that removal measure. We show that this problem is solvable for every compactly supported Borel probability measure.
	\end{abstract}
	
	\maketitle

	\section{Introduction and Main Results}
	
	We study a version of $n$-particle branching Brownian motion with drift in which each branching event is accompanied by the removal of a particle from a most populated spatial cell. This mechanism trims the population while keeping its size fixed, as in the branching random walk model of~\cite{AMW}. We consider the particle system together with its macroscopic density evolution, described by a parabolic free obstacle problem. Mass is removed from the top of the density profile, whose height serves as an upper obstacle. This height is not prescribed but is determined together with the density under the unit mass constraint; hence the term `free'.
	
	The model can be interpreted as describing a population of mobile organisms in a habitat with limited resources. As in earlier models of branching and selection~\cite{BDMM07}, this constraint motivates the assumption of a fixed population size. In our model, each birth is balanced by a death. Individuals reproduce at their current location and move during their lifetimes: Diffusion describes random dispersal,
	drift represents directed movement or transport by the environment. Models that combine movement, reproduction and competition have been studied in population dynamics~\cite{HL04}. Our removal rule represents competition for local resources: Individuals die in locations where crowding is greatest. Mortality also depends on local crowding in the model of~\cite{FM04}, where individuals are motionless and dispersal occurs at birth.
	
	Examples in~\cite{AMW} exhibit both absolutely continuous and atomic removal measures. Some involve discontinuous drifts and lie outside the scope of the well-posedness and convergence results established there. They motivate two questions that we pursue here: whether the theory extends to less regular drifts, and how rich a class of removal measures can be obtained by choosing different drifts.
	
	\subsection*{The Particle System}
	We consider $n$ particles on the real line. Between branching and selection events, the particles evolve independently as Brownian motions with drift~$b$. Each particle branches at rate one, producing an offspring at its current position. The exponential branching clocks of rate~1 are independent of one another and of the Brownian motions. Space is partitioned into cells of length~$\eps$. Each branching event is accompanied by the removal of a particle chosen uniformly from those occupying the most populated cells, so that the population size remains equal to $n$. The removal choice is based on the configuration immediately before branching and, conditional on this configuration, is independent of the choice of the branching particle.
	
	Equivalently, we may transfer the particle selected for removal to the position of the branching particle. Giving the offspring the label of the removed particle allows us to describe the system using a fixed set of $n$ labels. After each event, the diffusion motion resumes from the updated	positions, using the same Brownian motions for the fixed labels, and the	mechanism repeats.
	
	We now give a precise construction. Fix $n\in\N$ and $\eps>0$, and write
	$[n]=\{1,\dots,n\}$. Let $\mu_0$ be a probability measure on $\R$ with density $u_0\in\CC_0(\R,\R)$. Here~$\CC_0(\R,\R)$ denotes the space of continuous functions on~$\R$ that vanish at infinity. We assume that $b:\R\to\R$ is bounded and Borel measurable.	Define the cells
	$$
	C_k^\eps=[\eps k,\eps(k+1)),\qquad k\in\Z.
	$$
	Now consider a configuration $x=(x_1,\dots,x_n)\in\R^n$. Let
	$$
	\kappa_k^\eps(x)=|\{i\in[n]:x_i\in C_k^\eps\}|
	$$
	be the number of particles in cell $k$, and let
	$$
	A^\eps(x)
	=\operatorname{argmax}_{k\in\Z}\kappa_k^\eps(x)
	=\{k\in\Z:\kappa_k^\eps(x)\ge\kappa_{k'}^\eps(x)
	\text{ for all }k'\in\Z\}
	$$
	be the set of indices of the most populated cells. For any $k\in A^\eps(x)$,
	the total number of particles eligible for removal is
	$\kappa_k^\eps(x)|A^\eps(x)|$ so that uniform selection among these particles assigns each of them removal probability $1/(\kappa_k^\eps(x)|A^\eps(x)|)$.
	
	To describe the branching and selection events, we introduce a Poisson random measure~$\pi^n$ with intensity measure
	$$
	\hat\pi^n=\lambda_{(0,\iy)}\otimes\zeta_{[n]}
	\otimes\zeta_{[n]}\otimes\zeta_\Z\otimes\lambda_{[0,1]}
	$$
	on $(0,\iy)\times[n]\times[n]\times\Z\times[0,1]$.
	Here $\lambda_D$ and $\zeta_D$ denote Lebesgue measure and counting measure
	restricted to $D$, respectively. The marks $i$ and~$j$ will label the particle
	selected for removal and the branching particle, respectively; $k$ will be the cell containing particle $i$; and $\theta$ will be an auxiliary mark used to implement the removal probabilities. We set
	$$
	g^{n,\eps}(x;i,k,\theta)
	=\1_{\{\theta\kappa_k^\eps(x)|A^\eps(x)|\le1,\,
		x_i\in C_k^\eps,\,k\in A^\eps(x)\}},
	$$
	for $x\in\R^n$, $i\in[n]$, $k\in\Z$ and $\theta\in[0,1]$.		
	For a fixed configuration $x$, integrating $g^{n,\eps}$ over $\theta$
	gives rate $1/(\kappa_k^\eps(x)|A^\eps(x)|)$ for each candidate pair $(i,j)$
	with $x_i\in C_k^\eps$ and $k\in A^\eps(x)$. Summing over the eligible
	removal labels and their cells gives total rate~$1$ for each branching
	label~$j$. Thus the removal label is uniform over the eligible particles,
	and the Poisson construction realizes the mechanism described above.
	
	We now define the dynamics of the particles. Let $X^n(0)\sim\mu_0^{\otimes n}$ be the initial configuration, independent of~$\pi_n$, and let $W^n=(W_i^n)_{i\in[n]}$ be a family of $n$ independent standard Brownian motions, independent of $X^n(0)$ and~$\pi_n$.
	Denoting by $X_i^{n,\eps}(t)$ the position of particle $i \in [n]$ at time $t \ge 0$ and writing $X^{n,\eps}(t)=(X_i^{n,\eps}(t))_{i\in[n]}$,
	we define the stochastic process $X^{n,\eps;b,u_0}=X^{n,\eps}$ by
	\begin{align}
		\label{eqn:dfnX}
		X^{n,\eps}(t) ={}& X^n(0) + W^n(t) + \int_0^t (b(X_i^{n,\eps}(s)))_{i \in [n]} \D s\notag\\
		{}& +\int_{(0,t] \times [n] \times [n] \times \Z \times [0,1]} e_i(X_j^{n,\eps}(s-)-X_i^{n,\eps}(s-)) g^{n,\eps}(X^{n,\eps}(s-);i,k,\theta) \D\pi^n(s,i,j,k,\theta), \qquad t\ge 0.
	\end{align}
	Here $e_i$ is the $i$-th standard basis vector of $\R^n$. Existence and uniqueness
	of a global strong solution are established in Lemma~\ref{lem:SDEsol} below.
	The resulting process has paths in $\calD([0,\iy),\R^n)$.
	Equivalently, we have
	\begin{align*}
		X_i^{n,\eps}(t)
		={}&X_i^n(0)+W_i^n(t)+\int_0^t b(X_i^{n,\eps}(s))\D s\notag\\
		&+\int_{(0,t]\times[n]\times\Z\times[0,1]} (X_j^{n,\eps}(s-)-X_i^{n,\eps}(s-))g^{n,\eps}(X^{n,\eps}(s-);i,k,\theta)\D\pi_i^n(s,j,k,\theta), \qquad t\ge 0,\, i \in [n].
	\end{align*}
	Here $\pi_i^n$ is the restriction of $\pi^n$ to marks whose removal label
	is $i$. Each accepted event therefore replaces the position of particle $i$
	by that of particle $j$, leaving the other particle positions unchanged.
	The normalization in~$g^{n,\eps}$ ensures that branching and selection events
	occur at total rate $n$, independently of the configuration:
	\begin{equation}
		\label{eqn:unitInt}
		\int_{(0,t]\times[n]\times[n]\times\Z\times[0,1]}
		g^{n,\eps}(X^{n,\eps}(s-);i,k,\theta)
		\hat\pi^n(\D s,\D i,\D j,\D k,\D \theta)=nt,\qquad t\ge0.
	\end{equation}
	The case $i=j$ is allowed. It corresponds to replacing the branching particle by its offspring and leaves the configuration unchanged. More generally, an accepted event produces a jump only when the two positions differ; the total rate $n$ counts all accepted events.
	
	\subsection*{Empirical and removal measures}
	We study the limit as $n\to\iy$, with $\eps=\eps_n\searrow0$ along a suitable
	sequence. Let $\MM_1$ denote the Polish space of Borel probability measures
	on $\R$, equipped with the topology of weak convergence. We introduce the
	empirical measure process
	$$
	\mu^{n,\eps;b,u_0}= \mu^{n,\eps}: [0,\infty) \to \MM_1, \quad t \mapsto \mu_t^{n,\eps}= \frac{1}{n}\sum_{i \in [n]} \delta_{X_i^{n,\eps;b,u_0}(t)}
	$$
	associated with $X^{n,\eps}$. Recalling $X^{n,\eps}\in\calD([0,\iy),\R^n)$, we have $\mu^{n,\eps}\in\calD([0,\iy),\MM_1)$ a.s. We equip this path space with
	the usual Skorokhod $J_1$ topology and write time arguments as subscripts.
	For a path $\mu\in\calC([0,\iy),\MM_1)$, we call a measurable function
	$u:[0,\iy)\times\R\to\R$ a density of $\mu$ if
	$$
	\mu_t(\D x)=u(t,x)\D x,\qquad t\ge0,
	$$
	as an equality of measures.
	
	To identify the limiting selection rule, we also keep track of where removals occur. Selection from the most populated cells should become removal at maxima of the limiting density. We therefore introduce a second component, which records the space-time locations of removal, normalized by the population size. At a time of removal {$t$}, its spatial distribution is uniform over the particles occupying the most populated cells immediately before $t$. Namely, define the Borel measure $\beta^{n,\eps;b,u_0}=\beta^{n,\eps}$ on
	$[0,\iy)\times\R$ by
	$$
	\beta^{n,\eps}(D) = \frac{1}{n} \int_{(0,\infty) \times [n] \times [n] \times \Z \times [0,1]} \1_{\{(s,X_i^{n,\eps}(s-)) \in D\}} g^{n,\eps}(X^{n,\eps}(s-);i,k,\theta) \D\hat{\pi}^n(s,i,j,k,\theta),
	$$
	for Borel sets $D \in \BB([0,\infty) \times \R)$.
	This is the predictable compensator of the normalized removal measure. Its randomness
	comes from the configuration in the integrand. Let $\MM_{\nearrow1}$ be the
	set of Borel measures $\beta$ on $[0,\iy)\times\R$ satisfying
	$\beta([0,T]\times\R)=T$ for every $T\ge0$. Equivalently, these measures have
	Lebesgue measure as their time marginal. We endow this set with the topology
	of weak convergence on compact time horizons, i.e., the initial topology under
	the family $(\beta\mapsto\beta\abs_{[0,T]\times\R})_{T\ge0}$ of measure
	restrictions. The resulting space is Polish. By~\eqref{eqn:unitInt},
	we have $\beta^{n,\eps}\in\MM_{\nearrow1}$. For each $\beta\in\MM_{\nearrow1}$, the
	disintegration theorem gives a probability kernel $(\beta_t)_{t\ge0}$ such that
	$$
	\beta(\D t,\D x)=\beta_t(\D x)\D t.
	$$
	The kernel is determined only for a.e.~$t \ge 0$. When $\beta$ is random, we choose
	a version measurable also in the underlying outcome, and the identity holds
	a.s.; see, e.g., \cite[Theorem~1.25]{Kal17}.

	\subsection*{The free obstacle problem}
	The macroscopic evolution equation combines diffusion, drift, branching and removal. To describe its structure as an obstacle problem, write $h(t)=\sup_{x\in\R}u(t,x)$ for all $t \ge 0$. The density lies below the upper obstacle $h(t)$, and removal acts only where this level is attained. Below the obstacle, the density therefore satisfies, in the sense of distributions,
	$$
	\dot u=\tfrac12u''-(bu)'+u
	\qquad\text{on }\{u<h\},
	$$
	with $\dot u=\pl_tu$ and $u'=\pl_xu$. The height $h(t)$ is not prescribed: It is determined together with the density under the unit mass constraint $\int_\R u(t,x)\D x=1$. Following the terminology of~\cite{AMW}, we call this a parabolic free obstacle problem. We now make this description precise through a weak formulation.
	
	First, to express the selection rule, for a function $u \in \CC([0,\infty) \times \R,\R)$ such that $u\abs_{[0,T]\times\R}$ is bounded for every~$T<\iy$, set
	$$
	\BB^u=\{\beta\in\MM_{\nearrow1}:
	\beta(\{(t,x)\in[0,\iy)\times\R:
	u(t,x)\lne\y u(t,\cdot\,)\y_\infty\})=0\}.
	$$
	For a nonnegative density $u$, this condition means that $\beta_t$ is supported on the spatial maximizers of $u(t,\cdot)$ for a.e.~$t$.
	For $u$ as before and $\beta\in\MM_{\nearrow1}$, we call $(u,\beta)$ a weak solution of
	\begin{equation}
		\label{eqn:PDE}
		\dot{u} = \frac{1}{2} u'' - (bu)' + u - \beta, \qquad u(0, \cdot\,) = u_0, \qquad \beta \in \BB^u,
	\end{equation}
	if we have~$\beta\in\BB^u$ and the identity
	\begin{align*}
		\int_\R \varphi(x) (u(t,x) - u_0(x)) \D x
		={}& \int_0^t \int_\R \left(\frac{1}{2} \varphi''(x) + \varphi'(x) b(x) + \varphi(x) \right) u(s,x) \D x \D s
		- \int_{[0,t]\times \R} \varphi(x) \D \beta(s,x)
	\end{align*}
	holds for every $t\ge0$ and every test function $\varphi\in\CC_c^\infty(\R,\R)$.
	In this equation, the diffusion term is $\tfrac12u''$, transport is described by $-(bu)'$, branching contributes $u$, and $\beta$ describes removal.
	
	We work in the following class of continuous functions, whose members vanish at spatial infinity uniformly on bounded time intervals and have uniformly vanishing $L^2$-tails on those intervals. These conditions will be used in the regularization and uniqueness arguments:
	\begin{align*}
		\UU ={}& \left\{u \in \CC([0,\infty) \times \R,\R): \lim_{R \to \infty} \sup_{t \in [0,T]} \int_{[-R,R]^c} \abs u(t,x)\abs^2 \D x = 0 = \lim_{R \to \infty} \sup_{t \in [0,T]}\sup_{x \in [-R,R]^c} \abs u(t,x)\abs \text{ for all } T \ge 0\right\}.
	\end{align*}

	\subsection*{Stationary examples and motivating questions}
	
	To illustrate the role of the drift, recall two stationary examples from~\cite{AMW} (rescaled to match the diffusion coefficient in~\eqref{eqn:PDE}). For the smooth drift
	$$
	b(x)=-\sqrt{2}\tanh(\sqrt{2}x),
	$$
	set
	$$
	w=\frac{\log(1+\sqrt{2})}{\sqrt{2}},
	\qquad h=\frac{1}{2(w+1)}.
	$$
	A stationary solution is given by
	$$
	u(t,x)=
	\begin{cases}
		h, & |x|\le w,\\
		\displaystyle 2h\frac{\sinh(\sqrt{2}|x|)}{\cosh^2(\sqrt{2}|x|)},
		& |x|>w,
	\end{cases}
	$$
	and
	$$
	\beta_t(\D x)=h\left(1+\frac{2}{\cosh^2(\sqrt{2}x)}\right)
	\1_{[-w,w]}(x)\D x.
	$$
	The density has a flat top on $[-w,w]$, and the removal measure is absolutely continuous and supported on it.
	
	For the discontinuous inward drift
	$$
	b(x)=-\sqrt{2}\operatorname{sign}(x),
	$$
	a stationary solution is given by
	$$
	u(t,x)=\frac{\sqrt{2}}{4}
	(1+\sqrt{2}|x|)\e^{-\sqrt{2}|x|},
	\qquad \beta_t=\delta_0.
	$$
	The density has a unique maximizer at $0$, where the entire removal measure is concentrated. This example falls outside the scope of the well-posedness and convergence results in~\cite{AMW}, which require a Lipschitz drift.
	
	These examples led to the questions posed in~\cite[Problems 1.8 and 1.9]{AMW}: whether well-posedness and convergence theory extends beyond Lipschitz drift, and how the regularity of the drift constrains the maximizer sets and the removal measure. In the stationary setting, we approach the latter question through an inverse problem, asking which removal measures can be realized by a suitable choice of drift.
	
	\subsection*{Main Results}
	Under the standing assumptions on $b$ and $u_0$, the following theorem establishes existence and uniqueness for~\eqref{eqn:PDE} and identifies its solution as the deterministic hydrodynamic limit of the particle system. Convergence holds along a suitably chosen sequence of cell sizes $\eps_n\searrow0$. The path-space topology and the topology for the removal measures are those specified before.
	
	\begin{thm}
		\label{thm:1}
		There exists a null sequence~$(\eps_n)_{n \in \N} \subseteq (0,\infty)$ such that~$(\mu^{n,\eps_n},\beta^{n,\eps_n})_{n \in \N}\subseteq \DD([0,\infty),\MM_1) \times \MM_{\nearrow 1}$ converges in probability to a deterministic limit~$(\mu,\beta) \in \CC([0,\infty),\MM_1) \times \MM_{\nearrow 1}$ such that~$\mu$ has density~$u \in \UU$ and~$(u,\beta)$ is the unique solution of~\eqref{eqn:PDE} in~$\UU \times \MM_{\nearrow 1}$.
	\end{thm}
	
	For the inverse problem, we are given a target stationary removal measure $\beta_0$ and seek a drift and a stationary density that realize it. Assuming $\beta_0$ has compact support, we choose a density that is constant on an interval containing its support. In the interval, the diffusion term vanishes, and the stationary equation relates $\beta_0$ to the drift. We construct the drift using the distribution function of $\beta_0$, as specified in the following theorem.
	
	\begin{thm}
		\label{th2}
		Let~$\beta_0 \in \MM_1$ and~$x_-,x_+ \in \R$ with~$\beta_0([x_-,x_+])=1$. Setting $c=(\sqrt{2}+x_+-x_-)^{-1}$, we consider the probability density~$u_0 \in \CC_0(\R,\R)$ given by $$u_0:\R \to \R, \quad x \mapsto c \e^{-\sqrt{2} \op{dist}(x,[x_-,x_+])},$$
		and the càdlàg function~$$b:\R \to \R, \quad x \mapsto
		\begin{cases}
			\sqrt{2}, \quad& x \lne x_-,\\
			\left(x-\frac{x_++x_-}{2}\right)+\frac{1-2\beta_0((-\infty,x])}{2c}, \quad& x \in [x_-,x_+),\\
			-\sqrt{2}, \quad& x \ge x_+.
		\end{cases}$$
		Finally, let~$\mu \in \CC([0,\infty),\MM_1)$ have density~$u:[0,\infty) \times \R\to\R$, $(t,x)\mapsto u_0(x)$ and define $\beta=\lambda_{[0,\infty)}\otimes \beta_0 \in \MM_{\nearrow 1}$.
		Then there exists a null sequence~$(\eps_n)_{n \in \N} \subseteq (0,\infty)$ such that~$(\mu^{n,\eps_n;b,u_0},\beta^{n,\eps_n;b,u_0})_{n \in \N}\subseteq \DD([0,\infty),\MM_1) \times \MM_{\nearrow 1}$ converges in probability to~$(\mu,\beta)$ and~$(u,\beta)$ is the unique solution of~\eqref{eqn:PDE} in~$\UU \times \MM_{\nearrow 1}$.
	\end{thm}
	
	The constructed drift is bounded and c\`adl\`ag, hence satisfies our regularity assumption. Choosing $\beta_0$ to be the Cantor probability measure, for example, produces a stationary solution with a purely singular continuous removal measure. On~$[0,1]$, the corresponding drift is an affine function minus a multiple of the Cantor function. Stationarity here refers to the PDE solution: $u(t,\cdot)=u_0$ and $\beta_t=\beta_0$ for a.e.~$t\ge0$.
	
	Compact support is necessary within the stationary solution class considered
	here. A continuous probability density vanishing at infinity has a
	compact set of spatial maximizers, and the stationary removal measure is
	supported on that set. The inverse construction is not unique: for the same prescribed measure $\beta_0$, different choices of the containing interval $[x_-,x_+]$ yield
	different stationary densities and drifts. This is consistent with forward uniqueness in Theorem~\ref{thm:1},
	where $b$ and $u_0$ are fixed.
	
	\subsection*{Related work}
	
	Branching particle systems with selection provide models for populations subject to reproduction and competition. In models motivated by natural selection, particle positions represent fitness, and selection favors individuals with higher fitness; see, e.g., \cite{BDMM06}. A related model is the $N$-particle branching Brownian motion ($N$-BBM), in which each particle branches at rate one and the leftmost particle is removed upon each branching event, keeping the population size fixed. For suitable initial conditions, the fluctuations of the particle cloud on time scales of order $(\log N)^3$ were studied in~\cite{Mai16}.
	
	The hydrodynamic limit of the $N$-BBM is described by a free boundary problem, where the boundary represents the left edge of the population; see~\cite{DFPS19,BBP19}. The stationary empirical profile of the $N$-BBM, viewed from its leftmost particle, converges as $N\to\iy$ to the profile of the minimal travelling wave~\cite{BT24}. For suitable initial conditions, convergence to travelling waves and precise large time asymptotics for the front position of the free boundary problem were established in~\cite{BPT25}. Related connections between selection and free boundary problems arise for nonlocal branching models~\cite{DR11} and, in higher dimensions, for the Brownian bees model, where the particle farthest from the origin is removed upon each branching event~\cite{BBNP22}. Selection can also be encoded through a removal measure and a complementarity condition, without explicitly introducing a moving boundary. Such a weak formulation and its application to hydrodynamic limits were developed in~\cite{Ata25}.
	
	Spatial population models with mortality rates depending on local competition were studied in~\cite{FM04}. Different scalings yield convergence to a deterministic nonlinear integro-differential equation or to a nonlinear superprocess.
	
	In~\cite{AMW}, branching random walks in arbitrary spatial dimension are regulated by keeping the population size fixed and removing particles from maximally populated sites. Our Brownian particle model leads to the same type of limiting free obstacle problem. The main technical contribution concerns the drift: We allow very general bounded drifts, whereas the well-posedness and convergence results of~\cite{AMW} require Lipschitz continuity. Even in dimension one, admitting such drifts requires substantial changes to the regularity estimates and uniqueness argument.

	\subsection*{Outline}
	In Section~\ref{sec:prelim}, we discuss some preliminaries and reduce Theorem~\ref{th2} to Theorem~\ref{thm:1}. Section~\ref{sec:tight} deals with tightness while Section~\ref{sec:subseq} establishes properties of arbitrary subsequential distributional limits. In Section~\ref{sec:unique}, we prove that~\eqref{eqn:PDE} has at most one solution in~$\UU \times \MM_{\nearrow 1}$. Finally, in Section~\ref{sec:final}, the complementarity condition is verified and the proof of Theorem~\ref{thm:1} is completed.

	\section{Preliminaries and Proof of Theorem~\ref{th2}}
	\label{sec:prelim}
	We start by showing the elementary fact that the total number of branching-removal events of the~$n$ particles on any bounded time interval~$I$ is Poisson distributed with parameter~$n\lambda(I)$. Note that the number of jumps in that interval is bounded by the number of such events; in particular, removing the branching particle itself does not cause a jump.
	
	\begin{lem}
		\label{lem:jumpPoi}
		Let~$\eps \gne 0$ and let~$Y$ be a predictable $\R^n$-valued process and~$I \subseteq (0,\infty)$ a bounded interval. Then $$\int_{I \times [n] \times [n] \times \Z \times [0,1]} g^{n,\eps}(Y(s);i,k,\theta) \D\pi^n(s,i,j,k,\theta)$$ is Poisson distributed with parameter~$n\lambda(I)$.
	\end{lem}
	\begin{proof}
		Noting	$$\int_{I \times [n] \times [n] \times \Z \times [0,1]} \1_{\{\theta m_{k'}^\eps(Y(s)) \abs A^\eps(Y(s))\abs \le 1,\, Y_{i'}(s) \in C_{k'}^\eps,\, k' \in A^\eps(Y(s))\}} \D\hat{\pi}^n(s,i,j,k,\theta) = n\lambda(I) \quad \text{a.s.},$$
		the classic formula for the Laplace transform of a Poisson integral with predictable integrand yields
		\begin{align*}
			&\E \exp\left(-z \int_{I \times [n] \times [n] \times \Z \times [0,1]} g^{n,\eps}(Y(s);i,k,\theta) \D\pi^n(s,i,j,k,\theta)\right)\\
			={}& \E \exp\left(- \int_{I \times [n] \times [n] \times \Z \times [0,1]} (1- \e^{-zg^{n,\eps}(Y(s);i,k,\theta)}) \D\hat{\pi}^n(s,i,j,k,\theta)\right)\\
			={}& \E \exp\left(- (1- \e^{-z}) \int_{I \times [n] \times [n] \times \Z \times [0,1]} \1_{\{\theta m_{k'}^\eps(Y(s)) \abs A^\eps(Y(s))\abs \le 1,\, Y_{i'}(s) \in C_{k'}^\eps,\, k' \in A^\eps(Y(s))\}} \D\hat{\pi}^n(s,i,j,k,\theta)\right)\\
			={}& \exp(- (1- \e^{-z}) n\lambda(I)), \qquad z \ge 0,
		\end{align*}
		which is precisely the Laplace transform of a Poisson distribution with parameter~$n\lambda(I)$.
	\end{proof}
	
	It is mathematical folklore that, due to boundedness of~$b$, the jump diffusion~$X^{n,\eps}$ is well-defined by~\eqref{eqn:dfnX} for every~$\eps \gne 0$. In lack of a proper reference and for the sake of completeness, we briefly sketch the proof.
	\begin{lem}
		\label{lem:SDEsol}
		For every~$\eps \gne 0$, equation~\eqref{eqn:dfnX} has a unique global strong solution.
	\end{lem}
	\begin{proof}
		Let~$\tilde{x} \in \R^n$. According to a classic result~\cite{Ver80}, the SDE $$\tilde{X}^n(t) = \tilde{x} + W^n(t) + \int_0^t \big(b(\tilde{X}_i^n(s))\big)_{i \in [n]} \D s, \quad t \ge 0,$$ has a unique global strong solution~$\tilde{X}^n$.
		As a consequence of Lemma~\ref{lem:jumpPoi}, the number of jumps induced by the Poisson integral w.r.t.~$\pi^n$ in~\eqref{eqn:dfnX} is a.s.~finite on any compact time interval. Thus, by applying the above findings piecewise, \eqref{eqn:dfnX} must have a unique global strong solution as well.
	\end{proof}
	
	We end this section by deducing Theorem~\ref{th2} from Theorem~\ref{thm:1}.
	\begin{proof}[Proof of Theorem~\ref{th2}]
		By Theorem~\ref{thm:1}, there is a null sequence~$(\eps_n)_{n \in \N} \subseteq (0,\infty)$ such that~$(\mu^{n,\eps_n;b,u_0},\beta^{n,\eps_n;b,u_0})_{n \in \N}$ converges in probability to some deterministic~$(\tilde{\mu},\tilde{\beta}) \in \CC([0,\infty),\MM_1) \times \MM_{\nearrow 1}$ satisfying the following properties: There exists a function~$\tilde{u} \in \UU$ with~$\D \tilde{\mu}_t(x) = \tilde{u}(t,x) \D x$ for all~$t \in [0,\infty)$ as equality of measures. Further, $(\tilde{u},\tilde{\beta})$ is the unique solution of~\eqref{eqn:PDE} in~$\UU \times \MM_{\nearrow 1}$. Noting~$(u,\beta) \in \UU \times \MM_{\nearrow 1}$, it remains to prove that~$(u,\beta)$ is a solution of~\eqref{eqn:PDE} as well. As a consequence of~$\beta_0([x_-,x_+])=1$, we have~$\beta \in \BB^u$. Now let~$\varphi \in \CC_c^\infty(\R,\R)$. Since~$u$ and~$\beta$ are stationary in time, it suffices to prove
		\begin{align}
			\label{eqn:ODE}
			\int_\R \left(\frac{1}{2} \varphi''(x) + \varphi'(x) b(x) + \varphi(x) \right) u_0(x) \D x
			= \int_{\R} \varphi \D \beta_0
			= \int_{[x_-,x_+]} \varphi \D \beta_0.
		\end{align}
		Inserting the definitions of~$b$ and~$u_0$, we get
		\begin{align*}
			\int_{-\infty}^{x_-} \left(\frac{1}{2} \varphi''(x) + \varphi'(x) b(x) + \varphi(x) \right) u_0(x) \D x
			={}& \int_{-\infty}^{x_-} \left(\frac{1}{2} \varphi''(x) + \sqrt{2} \varphi'(x) + \varphi(x) \right) c\e^{\sqrt{2}(x-x_-)} \D x\\
			={}& \left[\left(\frac{1}{2} \varphi'(x) + \frac{1}{\sqrt{2}} \varphi(x)\right) c\e^{\sqrt{2}(x-x_-)}\right]_{x=-\infty}^{x_-}\\
			={}& \frac{c}{2} \varphi'(x_-) + \frac{\sqrt{2}c}{2} \varphi(x_-)
		\end{align*}
		and similarly
		\begin{align*}
			\int_{x_+}^\infty \left(\frac{1}{2} \varphi''(x) + \varphi'(x) b(x) + \varphi(x) \right) u_0(x) \D x
			={}& -\frac{c}{2} \varphi'(x_+) + \frac{\sqrt{2}c}{2} \varphi(x_+).
		\end{align*}
		Using integration by parts for Riemann and for Stieltjes integrals, we further obtain
		\begin{align*}
			&\int_{x_-}^{x_+} \left(\frac{1}{2} \varphi''(x) + \varphi'(x) b(x) + \varphi(x) \right) u_0(x) \D x\\
			={}& \int_{x_-}^{x_+} c \left(\frac{1}{2} \varphi''(x) + \varphi'(x) \left(x-\frac{x_++x_-}{2}+\frac{1}{2c}\right) + \varphi(x)\right) \D x - \int_{x_-}^{x_+} \beta_0((-\infty,x]) \varphi'(x) \D x\\
			={}& \frac{c}{2}[\varphi'(x)]_{x=x_-}^{x_+} + c\left[\varphi(x) \left(x-\frac{x_++x_-}{2}+\frac{1}{2c}\right)\right]_{x=x_-}^{x_+}
			- [\varphi(x)\beta_0((-\infty,x])]_{x=x_-}^{x_+} + \int_{[x_-,x_+]} \varphi \D \beta_0\\
			={}& \frac{c}{2} \varphi'(x_+)- \frac{c}{2} \varphi'(x_-) + \frac{(x_+-x_-)c}{2} (\varphi(x_+)+\varphi(x_-)) - \frac{1}{2} (\varphi(x_+)+\varphi(x_-)) + \int_{\R} \varphi \D \beta_0.
		\end{align*}
		Summing the three equations and inserting the definition of~$c$ yields~\eqref{eqn:ODE}.
	\end{proof}

	\section{Tightness}
	\label{sec:tight}
	Let us now turn toward the proof of Theorem~\ref{thm:1}. We start with tightness. Some parts of our analysis are based on a natural embedding of the underlying particle system into the following pure branching particle system, without removals or transfers: For each~$i \in [n]$,
	we independently start in~$X^n_i(0)$ a binary branching process of unit branching rate, with particles moving (with independent paths after branching, as usual) according to the diffusion $$\D Y(t) = \D W(t) + b(Y(t)) \D t, \quad t \ge 0,$$ where~$W$ is a standard Brownian motion. Note that, by~\cite{Ver80}, the above SDE has a unique global strong solution. Now fix~$i \in [n]$. We denote the number of particles alive at time~$t \ge 0$ in the $i$-th of these branching processes by~$N_i^n(t)$. Further, we enumerate the particle trajectories of that branching process, understood as paths on the whole time axis~$[0,\infty)$ (including the paths of the ancestors), and denote the resulting family of stochastic processes in~$\CC([0,\infty),\R)$ by~$(Y_{i,l}^n)_{l \in \N}$. W.l.o.g., we assume that the processes are enumerated such that~$\{Y_{i,l}^n(t): l \in [N_i^n(t)]\}$ is the set of all particle positions of the $i$-th branching process at time~$t$ for every~$t \ge 0$. Interpreting the particle transfers in our original system as removals with simultaneous branching of another particle, there is a natural coupling between the two systems. For any given~$\eps \gne 0$, this coupling satisfies
	$$\{(t,X_i^{n,\eps}(t)): t \ge 0,\, i \in [n]\} \subseteq \{(t,Y^n_{i,l}(t)): t \ge 0,\, i \in [n],\, l \in [N_i^n(t)]\} \quad \text{a.s.}$$
	
	We are now ready to state and prove our tightness result. Note that the cell size does not have to converge to~$0$.
	
	\begin{prop}
		\label{prop:tight}
		The family~$(\mu^{n,\eps},\beta^{n,\eps})_{n \in \N, \eps \gne 0}$ is tight. Further, for every sequence~$(\eps_n)_{n \in \N} \subseteq (0,\infty)$, each subsequential distributional limit~$\mu$ of~$(\mu^{n,\eps_n})_{n \in \N}$ a.s.~satisfies $\mu \in \CC([0,\infty),\MM_1)$.
	\end{prop}
	\noindent Since, for every sequence~$(\eps_n)_{n \in \N} \subseteq (0,\infty)$, the strong law of large numbers implies $$\lim_{n \to \infty} \mu_0^{n,\eps_n} = \lim_{n \to \infty} \frac{1}{n} \sum_{i \in [n]} \delta_{X_i^n(0)} = \mu_0 \quad \text{a.s.},$$ the notation in the proposition is consistent in the sense that the subsequential distributional limit~$\mu$ a.s.~starts in~$\mu_0$.
	\begin{proof}[Proof of Proposition~\ref{prop:tight}]
		As a consequence of Prohorov's theorem, tightness of~$(\mu^{n,\eps},\beta^{n,\eps})_{n \in \N, \eps \gne 0}$ is equivalent to tightness of~$(\mu^{n,\eps_n},\beta^{n,\eps_n})_{n \in \N}$ for every sequence~$(\eps_n)_{n \in \N} \subseteq (0,\infty)$. So let~$(\eps_n)_{n \in \N} \subseteq (0,\infty)$. It suffices to separately establish tightness of~$(\beta^{n,\eps_n})_{n \in \N}$ and $\CC$-tightness of~$(\mu^{n,\eps_n})_{n \in \N}$. In view of the topological structures of~$\DD([0,\infty),\MM_1)$ and~$\MM_{\nearrow 1}$, it further suffices to establish these properties on any compact time horizon~$[0,T]$. To this end, we fix~$T \gne 0$ for the rest of the proof. To simplify notation, we will refrain from indicating restrictions explicitly, i.e., we will write~$\mu^{n,\eps_n}$ and~$\beta^{n,\eps_n}$ instead of~$\mu^{n,\eps_n}\abs_{[0,T]}$ and~$\beta^{n,\eps_n}\abs_{[0,T] \times \R}$.
		
		Let us start with tightness of~$(\beta^{n,\eps_n})_{n \in \N}$. To this end, let~$\eta \gne 0$. For each~$m \in \N$, choose~$r_m \gne 0$ with $$\p\left(\sup_{t \in [0,T]} \abs Y^1_{1,1}(t)\abs \gne r_m\right) \le \frac{1}{m \e^T 2^m}.$$ Using the many-to-one lemma~\cite{HR17}, the expected number of paths in the pure branching system which leave the ball of radius~$r_m$ around~$0$ before time~$T$ can be bounded by
		\begin{align}
			\label{eqn:Ubound}
			\sum_{i \in [n]} \E\left[\sum_{l \in [N_i^n(T)]} \1_{\{\sup_{t \in [0,T]} \abs Y^n_{i,l}(t)\abs \gne r_m\}}\right]
			\le{}& \e^T \sum_{i \in [n]} \E \1_{\{\sup_{t \in [0,T]} \abs Y^n_{i,1}(t)\abs \gne r_m\}}\notag\\
			={}& n\e^T \p\left(\sup_{t \in [0,T]} \abs Y^1_{1,1}(t)\abs \gne r_m\right)
			\le \frac{n}{m2^m}, \qquad m \in \N.
		\end{align}
		Denoting by~$\MM_T([0,T] \times \R)$ the space of all Borel measures on~$[0,T] \times \R$ with total mass~$T$, endowed with the topology of weak convergence, we set~$O_m:=\{x \in \R: \abs x\abs \gne r_m\}$ and $$K_m:= \left\{\beta \in \MM_T([0,T] \times \R): \beta([0,T] \times O_m) \le \frac{1}{m}\right\}, \quad m \in \N.$$ Let~$m_0 \in \N$ with~$\sum_{m=m_0}^\infty \frac{1}{2^m} \le \eta$. Then~$K:= \overline{\bigcap_{m=m_0}^\infty K_m}$ is compact by Prohorov's theorem. For the time being, fix~$m \in \N$.
		Under the coupling described above, $\beta^{n,\eps_n}$ is concentrated on~$$\{(t,X_i^{n,\eps_n}(t-)): t \in [0,T],\, i \in [n]\} \subseteq \{(t,Y^n_{i,l}(t-)): t \in [0,T],\, i \in [n],\, l \in [N_i^n(t)]\}.$$ Further, along any given trajectory of the pure branching system, a.s.~at most one removal of a particle of the original system can happen. Thus the total number of removals at locations in~$O_m$ is bounded by the number of paths of the pure branching system which enter~$O_m$. In other words, we have
		\begin{align*}
			&\int_{(0,T] \times [n] \times [n] \times \Z \times [0,1]} \1_{\{\abs X_i^{n,\eps_n}(s-)\abs \gne r_m\}} g^{n,\eps_n}(X^{n,\eps_n}(s-);i,k,\theta) \D\pi^n(s,i,j,k,\theta)\\
			\le{}& \sum_{i \in [n]} \sum_{l \in [N_i^n(T)]} \1_{\{\sup_{t \in [0,T]} \abs Y^n_{i,l}(t)\abs \gne r_m\}} \qquad \text{a.s.}
		\end{align*}
		Together with~\eqref{eqn:Ubound}, we obtain
		\begin{align*}
			\E[\beta^{n,\eps_n}([0,T] \times O_m)]
			={}& \E\left[\frac{1}{n}\int_{(0,T] \times [n] \times [n] \times \Z \times [0,1]} \1_{\{\abs X_i^{n,\eps_n}(s-)\abs \gne r_m\}} g^{n,\eps_n}(X^{n,\eps_n}(s-);i,k,\theta) \D\hat{\pi}^n(s,i,j,k,\theta)\right]\\
			={}& \E\left[\frac{1}{n}\int_{(0,T] \times [n] \times [n] \times \Z \times [0,1]} \1_{\{\abs X_i^{n,\eps_n}(s-)\abs \gne r_m\}} g^{n,\eps_n}(X^{n,\eps_n}(s-);i,k,\theta) \D\pi^n(s,i,j,k,\theta)\right]\\
			\le {}& \E\left[\frac{1}{n} \sum_{i \in [n]} \sum_{l \in [N_i^n(T)]} \1_{\{\sup_{t \in [0,T]} \abs Y^n_{i,l}(t)\abs \gne r_m\}}\right]
			\le \frac{1}{m2^m}.
		\end{align*}
		We deduce
		\begin{align*}
			\p(\beta^{n,\eps_n} \not \in K)
			\le{}& \sum_{m=m_0}^\infty \p\left(\beta^{n,\eps_n}([0,T] \times O_m) \gne \frac{1}{m}\right)
			\le \sum_{m=m_0}^\infty m \E[\beta^{n,\eps_n}([0,T] \times O_m)]
			\le \sum_{m=m_0}^\infty \frac{1}{2^m} \le \eta, \quad n \in \N,
		\end{align*}
		proving tightness of~$(\beta^{n,\eps_n})_{n \in \N}$.
		
		Let us now turn toward~$(\mu^{n,\eps_n})_{n \in \N}$. We will invoke a classic characterization of $\CC$-tightness. With the notation from above and $$K_m':=\left\{\nu \in \MM_1: \nu(O_m) \le \frac{1}{m}\right\}, \quad m \in \N,$$ the set~$K':= \overline{\bigcap_{m=m_0}^\infty K_m'}$ is compact by Prohorov's theorem. For any~$t \in [0,T]$ and~$m \in \N$, the number of particles of the original system which are in~$O_m$ at time~$t$ is bounded by the number of paths of the pure branching system which enter~$O_m$ at all: $$n\mu^{n,\eps_n}_t(O_m) \le \sum_{i \in [n]} \sum_{l \in [N_i^n(T)]} \1_{\{\sup_{t \in [0,T]} \abs Y^n_{i,l}(t)\abs \gne r_m\}}.$$ Together with~\eqref{eqn:Ubound}, we get
		\begin{align}
			\label{eqn:muTight1}
			\sup_{t \in [0,T]} \p(\mu_t^{n,\eps_n} \not \in K')
			\le{}& \sup_{t \in [0,T]} \sum_{m=m_0}^\infty \p\left(\mu^{n,\eps_n}_t(O_m) \gne \frac{1}{m}\right)
			\le \sup_{t \in [0,T]} \sum_{m=m_0}^\infty \frac{m}{n} \E[n\mu^{n,\eps_n}_t(O_m)]\notag\\
			\le{}& \sum_{m=m_0}^\infty \frac{m}{n} \E\left[\sum_{i \in [n]}  \sum_{l \in [N_i^n(T)]} \1_{\{\sup_{t \in [0,T]} \abs Y^n_{i,l}(t)\abs \gne r_m\}}\right]
			\le \sum_{m=m_0}^\infty \frac{1}{2^m} \le \eta, \qquad n \in \N.
		\end{align}
		Given~$\delta \gne 0$, letting~$d_{\op{LP}}$ be the Lévy-Prohorov metric on~$\MM_1$ and setting~$$\Delta_\delta:=\{(s,t) \in [0,T]^2: 0 \le t-s \le \delta\},$$ we introduce moduli of continuity defined by$$w_{\op{LP}}(\mu,\delta):=\sup_{(s,t) \in \Delta_\delta} d_{\op{LP}}(\mu_s,\mu_t), \quad \mu \in \DD([0,T],\MM_1),$$ as well as $$w(f,\delta):=\sup_{(s,t) \in \Delta_\delta} \abs f(s)-f(t)\abs, \quad f \in \DD([0,T],\R).$$ By continuity of~$Y_{1,1}^1$, there must be a~$\delta \in (0,\frac{\eta^3}{32(T+1)} \wedge 1)$ with~$\p(w(Y_{1,1}^1,\delta) \gne \eta) \le \frac{\eta^2}{4}$.
		Denoting the $\eta$-neighborhood of a measurable set~$D \in \BB(\R)$ by~$D^\eta$, we obtain
		\begin{align*}
			&\p(w_{\op{LP}}(\mu^{n,\eps_n},\delta) \gne \eta)\notag\\
			={}& \p(d_{\op{LP}}(\mu^{n,\eps_n}_t,\mu^{n,\eps_n}_s) \gne \eta \text{ for some } (s,t) \in \Delta_\delta)\notag\\
			\le{}& \p(\mu^{n,\eps_n}_s(D) \gne \mu^{n,\eps_n}_t(D^\eta) + \eta \text{ or } \mu^{n,\eps_n}_t(D) \gne \mu^{n,\eps_n}_s(D^\eta) + \eta \text{ for some } (s,t) \in \Delta_\delta \text{ and } D \in \BB(\R))\notag\\
			\le{}& \p(\abs\{i \in [n]: \abs X_i^{n,\eps_n}(s) - X_i^{n,\eps_n}(t)\abs \gne \eta\}\abs \gne \eta n \text{ for some } (s,t) \in \Delta_\delta)\notag\\
			\le{}& \p\left(\abs\{i \in [n]: \abs X_i^{n,\eps_n}(s) - X_i^{n,\eps_n}(t)\abs \gne \eta \text{ and } X_i^{n,\eps_n} \in \CC([s,t],\R)\}\abs \gne \frac{\eta n}{2} \text{ for some } (s,t) \in \Delta_\delta\right)\notag\\
			& \quad + \p\left(\abs\{i \in [n]: X_i^{n,\eps_n} \not\in \CC([s,t],\R)\}\abs \gne \frac{\eta n}{2} \text{ for some } (s,t) \in \Delta_\delta\right), \qquad n \in \N.
		\end{align*}
		On the one hand, the choice of~$\delta$ implies
		\begin{align*}
			&\p\left(\abs\{i \in [n]: \abs X_i^{n,\eps_n}(s) - X_i^{n,\eps_n}(t)\abs \gne \eta \text{ and } X_i^{n,\eps_n} \in \CC([s,t],\R)\}\abs \gne \frac{\eta n}{2} \text{ for some } (s,t) \in \Delta_\delta\right)\notag\\
			\le{}& \p\left(\abs\{i \in [n]: \abs X_i^{n,\eps_n}(s) - X_i^{n,\eps_n}(t)\abs \gne \eta \text{ and } X_i^{n,\eps_n} \in \CC([s,t],\R) \text{ for some } (s,t) \in \Delta_\delta\}\abs \gne \frac{\eta n}{2}\right)\notag\\
			\le{}& \frac{2}{\eta n} \E\abs\{i \in [n]: \abs X_i^{n,\eps_n}(s) - X_i^{n,\eps_n}(t)\abs \gne \eta \text{ and } X_i^{n,\eps_n} \in \CC([s,t],\R) \text{ for some } (s,t) \in \Delta_\delta\}\abs\notag\\
			={}& \frac{2}{\eta n} \sum_{i \in [n]} \p(\abs X_i^{n,\eps_n}(s) - X_i^{n,\eps_n}(t)\abs \gne \eta \text{ and } X_i^{n,\eps_n} \in \CC([s,t],\R) \text{ for some } (s,t) \in \Delta_\delta)\notag\\
			\le{}& \frac{2}{\eta n} \sum_{i \in [n]} \p(\abs Y_{1,1}^1(s) - Y_{1,1}^1(t)\abs \gne \eta \text{ for some } (s,t) \in \Delta_\delta)
			= \frac{2}{\eta} \p(w(Y_{1,1}^1,\delta) \gne \eta)
			\le \frac{\eta}{2}, \qquad n \in \N.
		\end{align*}
		On the other hand, bounding the number of particles which jump by the total number of jumps in the third step, applying Lemma~\ref{lem:jumpPoi} in the forth step and using the choice of~$\delta$ in the final line, we obtain
		\begin{align*}
			&\limsup_{n \to \infty} \p\left(\abs\{i \in [n]: X_i^{n,\eps_n} \not\in \CC([s,t],\R)\}\abs \gne \frac{\eta n}{2} \text{ for some } (s,t) \in \Delta_\delta\right)\notag\\
			\le{}& \limsup_{n \to \infty} \sum_{l \in \left[\left\lceil \frac{T}{\delta}\right\rceil\right]} \p\left(\abs\{i \in [n]: X_i^{n,\eps_n} \not\in \CC([(l-1)\delta,l\delta],\R)\}\abs \gne \frac{\eta n}{4}\right)\notag\\
			\le{}& \limsup_{n \to \infty} \frac{16}{\eta^2 n^2} \sum_{l \in \left[\left\lceil \frac{T}{\delta}\right\rceil\right]} \E\abs\{i \in [n]: X_i^{n,\eps_n} \not\in \CC([(l-1)\delta,l\delta],\R)\}\abs^2\notag\\
			\le{}& \limsup_{n \to \infty} \frac{16}{\eta^2 n^2} \sum_{l \in \left[\left\lceil \frac{T}{\delta}\right\rceil\right]} \E\left(\int_{((l-1)\delta,l\delta] \times [n] \times [n] \times \Z \times [0,1]} g^{n,\eps_n}(X^{n,\eps_n}(s-);i,k,\theta) \D\pi^n(s,i,j,k,\theta)\right)^2\notag\\
			={}& \limsup_{n \to \infty} \frac{16}{\eta^2 n^2} \sum_{l \in \left[\left\lceil \frac{T}{\delta}\right\rceil\right]} (n\delta + n^2 \delta^2)
			\le \frac{16(T+1)\delta}{\eta^2}
			\le \frac{\eta}{2}.
		\end{align*}
		Combining the three chains of inequalities, we get~$$\limsup_{n \to \infty} \p(w_{\op{LP}}(\mu^{n,\eps_n},\delta) \gne \eta) \le \eta.$$
		By a classic characterization (cf.~\cite[Corollary~VI.3.26]{EK86} and~\cite[Proposition~VI.3.26]{JS03}), this estimate together with~\eqref{eqn:muTight1} precisely gives the desired $\CC$-tightness of~$(\mu^{n,\eps_n})_{n \in \N}$.
	\end{proof}

	\section{Properties of Subsequential Hydrodynamic Limits}
	\label{sec:subseq}
	We start by showing that any subsequential distributional limit~$\mu$ of the particle distributions a.s.~has a density in~$\UU$. Further, we show that this density is dominated by the particle density corresponding to the pure branching system introduced in Section~\ref{sec:tight}. While we will only need the result for fixed cell size~$\eps \gne 0$, note that it is valid for any sequence~$(\eps_n)_{n \in \N} \subseteq (0,\infty)$.
	As usual, we denote the generic sample space encoding the underlying randomness by~$\Omega$.
	
	\begin{lem}
		\label{lem:density}
		Given~$(\eps_n)_{n \in \N} \subseteq (0,\infty)$, let~$\mu$ be a subsequential distributional limit of~$(\mu^{n,\eps_n})_{n \in \N}$.
		Then there exist a deterministic~$v \in \UU$, depending only on~$(u_0,b)$, and a product-measurable function $$u: [0,\infty) \times \R \times \Omega \to [0,\infty)$$ which a.s.~satisfies~$u \le v$ and~$\D \mu_t(x) = u(t,x) \D x$ for all~$t \ge 0$. Here, as usual, we suppress the third parameter of~$u$, which is governing its randomness.
		Further, the function~$v$ may be chosen with~$\int_\R v(t,x)\D x=\e^t$ for all~$t\ge0$ and with
		\begin{equation}
			\label{eqn:density-occupation}
			\mathbb E\left[\int_0^t\int_\R g(s,x)\D\mu_s^{n,\eps_n}(x)\D s\right] \le \int_0^t\int_\R g(s,x)v(s,x)\D x\D s, \quad n \in \N,\, t \ge 0,
		\end{equation}
		for every non-negative measurable function~$g:[0,\infty)\times\R\to[0,\infty)$.
	\end{lem}
	
	\begin{proof}
		Given~$Y(0) \sim \mu_0$, let~$Y=(Y(t))_{t \ge 0}$ be the unique global strong solution to~$$\begin{aligned}\D Y(t) = \D W^1(t) + b(Y(t)) \D t, \quad t \ge 0.\end{aligned}$$ Then there exists a~$\tilde{v} \in \CC_b([0,\infty) \times \R, \R)$ with~$\p(Y(t) \in D) = \int_D \tilde{v}(t,x) \D x$ for all~$t \ge 0$ and~$D \in \BB(\R)$. (In fact, the density~$\tilde{v}$ is the unique solution to the Fokker-Planck equation associated with the above SDE.)
		Let~$(Q_t)_{t \ge 0}$ be the transition semigroup of a standard Brownian motion and let~$T \ge 0$. It is easy to check from~\cite[Theorem~1]{QZ02} that there is a constant~$C_T \ge 0$ with~$\tilde{v}(t,\cdot) \le C_T Q_{2t}u_0$ for all~$t \in [0,T]$. Since~$[0,T] \to \CC_0(\R,\R)$, $t\mapsto Q_{2t} u_0$ is continuous, a compactness argument yields~$\lim_{R \to \infty} \sup_{t \in [0,T]} \sup_{x \in [-R,R]^c} \abs \tilde{v}(t,x)\abs = 0$. Moreover, we obtain
		\begin{align*}
			\sup_{t \in [0,T]} \int_{[-R,R]^c} \abs \tilde{v}(t,x)\abs^2 \D x
			\le{}& \y \tilde{v}\y_\infty \sup_{t \in [0,T]} \p(\abs Y(t)\abs \ge R)\\
			\le{}& \y \tilde{v}\y_\infty \p\left(\abs Y(0)\abs \ge \tfrac{R}{2}\right) + \y \tilde{v}\y_\infty \sup_{t \in [0,T]} \p\left(\abs W^1(t)\abs \ge \tfrac{R}{2} - \y b\y_\infty T\right)
			\xrightarrow{R \to \infty} 0.
		\end{align*}
		We deduce that~$\tilde{v}$ and hence also~$v: [0,\infty) \times \R \to \R$, $(t,x) \mapsto \e^t \tilde{v}(t,x)$ is an element of~$\UU$.
		
		Now denote the space of finite Borel measures on~$\R$, endowed with the topology of weak convergence, by~$\MM_f$. Using the notation for the pure branching system introduced in Section~\ref{sec:tight}, let~$\xi^n \in \DD([0,\infty),\MM_f)$ be its empirical measure process, defined by~$$\xi^n_t := \frac{1}{n} \sum_{i \in [n]} \sum_{l \in [N_i^n(t)]} \delta_{Y_{i,l}^n(t)}, \quad t \ge 0.$$
		By classic branching process theory, $(\xi^n)_{n \in \N}$ converges in probability to a deterministic limit~$\xi \in \CC([0,\infty),\MM_f)$, which satisfies~$\D \xi_t(x) = v(t,x) \D x$ for all~$t \ge 0$.
		Noting~$\mu^{n,\eps_n} \le \xi^n$ for every~$n \in \N$ under the coupling described in Section~\ref{sec:tight}, Strassen's theorem together with (subsequential) distributional convergence implies the existence of a coupling~$(\tilde{\mu},\tilde{\xi})$ of the limits~$\mu$ and~$\xi$ which satisfies~$\tilde{\mu} \le \tilde{\xi}$~a.s. Here the order is preserved under convergence to continuous paths. In particular, we a.s.~have~$\mu_t \le \xi_t$ and hence~$\mu_t \ll \lambda$ for all~$t \ge 0$. Thus~\cite[Theorem~V.58]{DM82} and the subsequent remark guarantee the existence of a product-measurable~$\check{u}: [0,\infty) \times \R \times \Omega \to [0,\infty)$ which a.s.~satisfies~$\D \mu_t(x) = \check{u}(t,x) \D x$ for all~$t \ge 0$. Recalling~$\mu_t \le \xi_t$ for all~$t \ge 0$ a.s., the random function~$u:=\check{u} \wedge v$ has the desired properties.
		
		Finally, let~$t \ge 0$ and~$n \in \N$. Since~$\tilde{v}(t,\cdot)$ is a probability density, we have~$\int_\R v(t,x)\D x=\e^t$. Moreover, the coupling~$\mu^{n,\eps_n}\le\xi^n$ as well as the first-moment formula for the pure branching system and Tonelli's theorem give
		\begin{align*}
			\mathbb E\left[\int_0^t\int_\R g(s,x)\D\mu_s^{n,\eps_n}(x)\D s\right]
			\le{}&\mathbb E\left[\int_0^t\int_\R g(s,x)\D\xi_s^n(x)\D s\right]
			= \int_0^t\int_\R g(s,x)v(s,x)\D x\D s	
		\end{align*}
		for every non-negative measurable function~$g:[0,\infty)\times\R\to[0,\infty)$.
	\end{proof}

	The next fundamental step is to establish a PDE, in weak form, for any subsequential distributional limit. We denote the partial derivative of a function~$\varphi \in \CC^{1,1}([0,\infty) \times \R,\R)$ w.r.t.~its first component (representing time) by~$\dot{\varphi}$ and the partial derivative w.r.t.~its second component (representing space) by~$\varphi'$, similar to the notation used in the formal PDE~\eqref{eqn:PDE}.
	
	\begin{lem}
		\label{lem:PDE}
		Given~$(\eps_n)_{n \in \N} \subseteq (0,\infty)$, let~$(\mu,\beta)$ be a subsequential distributional limit of~$(\mu^{n,\eps_n},\beta^{n,\eps_n})_{n \in \N}$. Then we a.s.~have
		\begin{align*}
			&\int_\R \varphi(t,x) \D \mu_t(x)\\
			={}& \int_\R \varphi(0,x) \D \mu_0(x)
			- \int_{[0,t]\times \R} \varphi(s,x) \D \beta(s,x)
			+ \int_0^t \int_\R \left(\frac{1}{2} \varphi''(s,x) + \varphi'(s,x) b(x) + \varphi(s,x) + \dot{\varphi}(s,x)\right) \D \mu_s(x) \D s
		\end{align*}
		simultaneously for all~$\varphi \in \CC_c^{1,2}([0,\infty) \times \R,\R)$ and~$t \ge 0$.
	\end{lem}
	
	\begin{proof}
		Let~$\varphi \in \CC_c^{1,2}([0,\infty) \times \R,\R)$ and~$t \ge 0$ be fixed. Further, let~$n \in \N$ and let~$c \in [0,\infty)$ be a bound on~$\abs \varphi\abs$ and~$\abs \varphi' \abs$. We define $$M^{n,1}(t) := \frac{1}{n} \sum_{i \in [n]} \int_0^t \varphi'(s,X_i^{n,\eps_n}(s)) \D W^n_i(s)$$ and
		\begin{align*}
			M^{n,2}(t) :={}& \frac{1}{n} \int_{(0,t] \times [n] \times [n] \times \Z \times [0,1]} (\varphi(s,X_j^{n,\eps_n}(s-)) - \varphi(s,X_i^{n,\eps_n}(s-))) g^{n,\eps_n}(X^{n,\eps_n}(s-);i,k,\theta) \D\pi^n(s,i,j,k,\theta)\\
			& - \frac{1}{n} \int_{(0,t] \times [n] \times [n] \times \Z \times [0,1]} (\varphi(s,X_j^{n,\eps_n}(s-)) - \varphi(s,X_i^{n,\eps_n}(s-))) g^{n,\eps_n}(X^{n,\eps_n}(s-);i,k,\theta) \D\hat{\pi}^n(s,i,j,k,\theta).
		\end{align*}
		Observing
		\begin{align*}
			&\frac{1}{n} \int_{(0,t] \times [n] \times [n] \times \Z \times [0,1]} \varphi(s,X_j^{n,\eps_n}(s-)) g^{n,\eps_n}(X^{n,\eps_n}(s-);i,k,\theta) \D\hat{\pi}^n(s,i,j,k,\theta)\\
			={}& \int_0^t \left(\frac{1}{n} \sum_{j \in [n]} \varphi(s,X_j^{n,\eps_n}(s-))\right) \left(\int_0^1 \sum_{i \in [n], k \in \Z} g^{n,\eps_n}(X^{n,\eps_n}(s-);i,k,\theta) \D \theta \right) \D s\\
			={}& \int_0^t \left(\frac{1}{n} \sum_{j \in [n]} \varphi(s,X_j^{n,\eps_n}(s))\right) 1 \D s\\
			={}& \int_0^t \int_\R \varphi(s,x) \D \mu_s^{n,\eps_n}(x) \D s
		\end{align*}
		and
		\begin{align*}
			\frac{1}{n} \int_{(0,t] \times [n] \times [n] \times \Z \times [0,1]} \varphi(s,X_i^{n,\eps_n}(s-)) g^{n,\eps_n}(X^{n,\eps_n}(s-);i,k,\theta) \D\hat{\pi}^n(s,i,j,k,\theta)
			=\int_{[0,t]\times \R} \varphi(s,x) \D \beta^{n,\eps_n}(s,x),
		\end{align*}
		Itô's formula for jump diffusions yields
		\begin{align}
			\label{eqn:Ito}
			&\int_\R \varphi(t,x) \D \mu^{n,\eps_n}_t(x)\notag\\
			={}& \frac{1}{n} \sum_{i \in [n]} \varphi(t,X_i^{n,\eps_n}(t))\notag\\
			={}& \frac{1}{n} \sum_{i \in [n]} \varphi(0,X_i^n(0))
			+ \frac{1}{n} \sum_{i \in [n]} \int_0^t \varphi'(s,X_i^{n,\eps_n}(s)) \D W^n_i(s)
			+ \frac{1}{n} \sum_{i \in [n]} \int_0^t \frac{1}{2} \varphi''(s,X_i^{n,\eps_n}(s)) \D s\notag\\
			&\quad + \frac{1}{n} \sum_{i \in [n]} \int_0^t \dot{\varphi}(s,X_i^{n,\eps_n}(s)) \D s + \frac{1}{n} \sum_{i \in [n]} \int_0^t \varphi'(s,X_i^{n,\eps_n}(s)) b(X_i^{n,\eps_n}(s)) \D s\notag\\
			&\quad +\frac{1}{n} \int_{(0,t] \times [n] \times [n] \times \Z \times [0,1]} (\varphi(s,X_j^{n,\eps_n}(s-)) - \varphi(s,X_i^{n,\eps_n}(s-))) g^{n,\eps_n}(X^{n,\eps_n}(s-);i,k,\theta) \D\pi^n(s,i,j,k,\theta)\notag\\
			={}& \int_\R \varphi(0,x) \D \mu_0^n(x) + \int_0^t \int_\R \left(\frac{1}{2} \varphi''(s,x) + \varphi(s,x) + \dot{\varphi}(s,x)\right) \D \mu_s^{n,\eps_n}(x) \D s\notag\\
			&\quad + \int_0^t \int_\R \varphi'(s,x) b(x) \D \mu_s^{n,\eps_n}(x) \D s -  \int_{[0,t]\times \R} \varphi(s,x) \D \beta^{n,\eps_n}(s,x) + M^{n,1}(t) + M^{n,2}(t).
		\end{align}
		Noting that the integrand in the Poisson integral defining~$M^{n,2}(t)$ is predictable, the (predictable) quadratic variations satisfy
		\begin{align*}
			\langle M^{n,1}\rangle(t)
			={}& \frac{1}{n^2} \sum_{i \in [n]} \int_0^t \abs \varphi'(s,X_i^{n,\eps_n}(s)) \abs^2 \D s
			\le \frac{c^2 t}{n} \xrightarrow{n \to \infty} 0
		\end{align*}
		and
		\begin{align*}
			\langle M^{n,2}\rangle(t)
			={}& \frac{1}{n^2} \int_{(0,t] \times [n] \times [n] \times \Z \times [0,1]} (\varphi(s,X_j^{n,\eps_n}(s-)) - \varphi(s,X_i^{n,\eps_n}(s-)))^2 g^{n,\eps_n}(X^{n,\eps_n}(s-);i,k,\theta) \D\hat{\pi}^n(s,i,j,k,\theta)\\
			\le{}& \frac{4 c^2}{n^2} \int_{(0,t] \times [n] \times [n] \times \Z \times [0,1]} g^{n,\eps_n}(X^{n,\eps_n}(s-);i,k,\theta) \D\hat{\pi}^n(s,i,j,k,\theta)
			= \frac{4 c^2 t}{n}
			\xrightarrow{n \to \infty} 0.
		\end{align*}
		Together with~$M^{n,1}(0)=M^{n,2}(0)=0$, this implies that~$M^{n,1}(t)$ and~$M^{n,2}(t)$ converge to~$0$ in probability as~$n \to \infty$. Consequently, Slutsky's theorem implies that~$(\mu,\beta,0,0)$ is a subsequential distributional limit of~$(\mu^{n,\eps_n},\beta^{n,\eps_n},M^{n,1}(t),M^{n,2}(t))_{n \in \N}$.
		Invoking Skorokhod's representation theorem, we assume a.s.~convergence.
		Since Proposition~\ref{prop:tight} guarantees~$\mu \in \CC([0,\infty),\MM_1)$ a.s., we a.s.~have convergence of~$(\mu^{n,\eps_n}_s)_{n \in \N}$ to~$\mu_s$ in~$\MM_1$ for all~$s \ge 0$ along the considered subsequence, by classic properties of Skorokhod spaces. All integrands in~\eqref{eqn:Ito} except the one involving~$b$ are bounded and continuous.
		For every bounded measurable function~$\tilde{b}:\R\to\R$, set
		$$I^n(\tilde{b}):=\int_0^t\int_\R\varphi'(s,x)\tilde{b}(x)\D\mu_s^{n,\eps_n}(x)\D s, \qquad I(\tilde{b}):=\int_0^t\int_\R \varphi'(s,x)\tilde{b}(x)\D\mu_s(x)\D s.$$ Assume for a moment that~$(I^n(b))_{n \in \N}$ converges to~$I(b)$ in probability along the considered subsequence. Then taking the subsequential limit in probability in~\eqref{eqn:Ito}, with the help of the dominated convergence theorem, yields the claimed formula for any fixed~$\varphi$ and~$t$. Since~$\CC_c^{1,2}([0,\infty) \times \R,\R)$ (w.r.t.~its usual topology) and~$[0,\infty)$ are separable, a continuity argument invoking the dominated convergence theorem guarantees the claimed uniformity of the a.s.-exceptional set.
		
		It remains to prove the subsequential convergence of~$(I^n(b))_{n \in \N}$ to~$I(b)$. Let~$u$ and~$v$ be as in Lemma~\ref{lem:density}, also satisfying the additional properties mentioned there.		
		Further, let~$\rho\in\CC_c^\infty(\R,[0,\infty))$ with~$\int_\R\rho(x)\D x=1$ and~$k \in \N$. Then~$b_k:=b*(k \rho(k \,\cdot\,))$ is continuous with~$\|b_k\|_\infty\le\|b\|_\infty$.
		By the subsequential convergence of~$\mu_s^{n,\eps_n}$ to~$\mu_s$ for every~$s \ge 0$ and the dominated convergence theorem, we have~$I^n(b_k) \to I(b_k)$ a.s.~as~$n \to \infty$ along the considered subsequence.
		Noting~$\int_0^t\int_\R v(s,x)\D x\D s=\e^t-1$ and~$\lim_{k \to \infty} b_k= b$ Lebesgue-a.e.~by Lebesgue's differentiation theorem, \eqref{eqn:density-occupation} and the dominated convergence theorem imply
		\begin{align*}
			\sup_{n \in \N} \E \abs I^n(b) - I^n(b_k)\abs
			\le{}& \sup_{n \in \N} \E\left[\int_0^t \int_\R c\abs b(x)-b_k(x)\abs \D\mu_s^{n,\eps_n}(x) \D s\right] \le \int_0^t\int_\R c\abs b(x)-b_k(x)\abs v(s,x)\D x\D s
			\xrightarrow{k \to \infty} 0.
		\end{align*}
		Invoking properties of~$u$, we get
		\begin{align*}
			\abs I(b_k) - I(b)\abs
			\le \int_0^t\int_\R c \abs b_k(x)-b(x)\abs u(s,x) \D x\D s
			\le \int_0^t\int_\R c \abs b_k(x)-b(x)\abs v(s,x) \D x\D s
			\xrightarrow{k \to \infty} 0 \quad \text{a.s.}
		\end{align*}
		Using the triangle inequality and Markov's inequality, the desired convergence in probability follows by first letting~$n\to\infty$ along the considered subsequence and then~$k\to\infty$.
	\end{proof}

	We end this section by essentially showing uniform equicontinuity of all possible subsequential particle densities.
	\begin{lem}
		\label{lem:cont}
		Let~$\beta \in \MM_{\nearrow 1}$ and let~$u: [0,\infty) \times \R \to \R$ be measurable with
		\begin{align*}
			&\int_\R \varphi(t,x) u(t,x) \D x\\
			={}& \int_\R \varphi(0,x) u_0(x) \D x
			- \int_{[0,t]\times \R} \varphi(s,x) \D \beta(s,x)\\
			&{}\quad+ \int_0^t \int_\R \left(\frac{1}{2} \varphi''(s,x) + \varphi'(s,x) b(x) + \varphi(s,x) + \dot{\varphi}(s,x)\right) u(s,x) \D x \D s, \qquad \varphi \in \CC_c^\infty([0,\infty) \times \R,\R) ,\, t \ge 0.
		\end{align*}
		Assume that, for every~$T \ge 0$, there is a constant~$C_T \in [0,\infty)$ with~$\y u\abs_{[0,T] \times \R} \y_\infty \le C_T$.
		Then there exists a function~$\tilde{u} \in \CC([0,\infty)\times \R,\R)$ with~$\tilde{u}(t,x)=u(t,x)$ for a.e.~$x \in \R$ and all~$t \ge 0$. Further, for every~$T \ge 0$, the restriction~$\tilde{u}\abs_{[0,T] \times \R}$ has a modulus of continuity, which depends only on~$(u_0,\y b\y_\infty,T,C_T)$, but not on~$\beta$ or~$u$ itself.
	\end{lem}
	\begin{proof}
		To simplify notation, we write~$u_s:= u(s,\cdot)$ for~$s \gne 0$. Let~$(Q_t)_{t \ge 0}$ be the transition semigroup of a standard Brownian motion. For the time being, we fix~$t \gne 0$ and~$\psi \in \CC_c^\infty(\R,\R)$. Let~$\varphi \in \CC^\infty_b([0,\infty) \times \R, \R)$ be a function satisfying~$\varphi(s,x):= \e^{t-s} Q_{t-s} \psi(x)$ for all~$(s,x) \in [0,t] \times \R$. Further, let~$R \gne 0$ and~$\chi_R \in \CC^\infty_c(\R,[0,1])$ with~$\chi_R\abs_{[-R,R]}=1$ as well as~$\y \chi_R'\y_\infty\le 1$ and~$\y \chi_R''\y_\infty \le 1$. Then $\varphi_R: [0,\infty) \times \R \to \R$, $(s,x) \mapsto \chi_R(x) \varphi(s,x)$ satisfies~$\varphi_R \in \CC_c^\infty([0,\infty) \times \R,\R)$. Inserting this test function into the given PDE in weak form and using the PDE~$\dot{\varphi}= -\frac{1}{2}\varphi''-\varphi$ on~$[0,t] \times \R$, we obtain
		\begin{align}
			\label{eqn:cutoff}
			&\int_\R \chi_R(x) \psi(x) u(t,x) \D x\notag\\
			={}& \int_\R \varphi_R(0,x) u_0(x) \D x
			- \int_{[0,t]\times \R} \varphi_R(s,x) \D \beta(s,x)\notag\\
			&\quad{}+ \int_0^t \int_\R \left(\frac{1}{2} \varphi_R''(s,x) + \varphi_R'(s,x) b(x) - \frac{1}{2}\chi_R(x)\varphi''(s,x)\right) u_s(x) \D x \D s \notag\\
			={}& \int_\R \chi_R(x) \varphi(0,x) u_0(x) \D x
			- \int_{[0,t]\times \R} \chi_R(x) \varphi(s,x) \D \beta(s,x)\notag\\
			&\quad{}+ \int_0^t \int_\R \left(\frac{1}{2} \chi_R''(x) \varphi(s,x) + \chi_R'(x) \varphi'(s,x) + \chi_R'(x) \varphi(s,x) b(x) + \chi_R(x) \varphi'(s,x) b(x)\right) u_s(x) \D x \D s.
		\end{align}
		Since~$(Q_t)_{t \ge 0}$ is a contraction semigroup on~$L^1(\R,\R)$, the functions $\varphi\abs_{[0,t] \times \R}$ and~$\varphi'\abs_{[0,t] \times \R}$ are integrable. Applying the dominated convergence theorem, we deduce
		$$\lim_{R \to \infty} \int_0^t \int_\R \left(\frac{1}{2} \chi_R''(x) \varphi(s,x) + \chi_R'(x) \varphi'(s,x) + \chi_R'(x) \varphi(s,x) b(x) \right) u_s(x) \D x \D s = 0.$$
		Since~$Q_t$ is self-adjoint (or using Fubini's theorem), we get
		\begin{align*}
			\int_\R \chi_R(x) \varphi(0,x) u_0(x) \D x
			={}& \int_\R \e^t Q_t \psi(x) \chi_R(x) u_0(x) \D x
			= \int_\R \e^t \psi(x) Q_t (\chi_R u_0)(x) \D x
			\xrightarrow{R \to \infty} \int_\R \e^t \psi(x) Q_t u_0(x) \D x,
		\end{align*}
		where the last step follows from the fact that the dominated convergence theorem implies $$\y Q_t((1-\chi_R)u_0)\y_1 \le \y (1-\chi_R)u_0\y_1 \xrightarrow{R \to \infty} 0.$$
		Together with integration by parts, a similar argument yields
		\begin{align*}
			\lim_{R \to \infty} \int_0^t \int_\R \chi_R(x) \varphi'(s,x) b(x) u_s(x) \D x \D s
			={}& \lim_{R \to \infty} \int_0^t  \int_\R \e^{t-s} Q_{t-s} \psi'(x) (\chi_R b u_s)(x) \D x \D s\\
			={}& \lim_{R \to \infty} \int_0^t \int_\R \e^{t-s} \psi'(x) Q_{t-s}(\chi_R b u_s)(x) \D x \D s\\
			={}& \int_0^t \int_\R \e^{t-s} \psi'(x) Q_{t-s} (bu_s)(x) \D x \D s\\
			={}& -\int_0^t \int_\R \e^{t-s} \psi(x) (Q_{t-s} (bu_s))'(x) \D x \D s.
		\end{align*}
		Consequently, taking the limit~$R \to \infty$ in~\eqref{eqn:cutoff}, we obtain
		\begin{align}
			\label{eqn:contIntermed}
			\int_\R \psi(x) u(t,x) \D x
			={}& \int_\R \e^t \psi(x) Q_t u_0(x) \D x
			- \int_{[0,t]\times \R} \varphi(s,x) \D \beta(s,x)
			-\int_0^t \int_\R \e^{t-s} \psi(x) (Q_{t-s} (bu_s))'(x) \D x \D s.
		\end{align}
		Now let~$(q_t)_{t \gne 0}$ be the family of transition densities of a standard Brownian motion.
		By Fubini's theorem, we have
		\begin{align*}
			\int_{[0,t]\times \R} \varphi(s,x) \D \beta(s,x)
			={}& \int_0^t \int_\R \int_\R \e^{t-s} q_{t-s}(x-y) \psi(y) \D y \D \beta_s(x) \D s
			= \int_0^t \int_\R \int_\R \e^{t-s} q_{t-s}(y-x) \D \beta_s(y) \psi(x)\D x \D s.
		\end{align*}
		Inserting this into~\eqref{eqn:contIntermed} and recalling that the resulting equation holds for all~$\psi \in \CC_c^\infty(\R,\R)$, we get~$u(t,\cdot) = \tilde{u}(t,\cdot)$ $\lambda$-a.e.~for $$\tilde{u}: [0,\infty) \times \R \to \R, \quad (t,x) \mapsto \e^t Q_t u_0(x)
		- \int_0^t \int_\R \e^{t-s} q_{t-s}(y-x) \D \beta_s(y) \D s
		-\int_0^t \e^{t-s} (Q_{t-s} (bu_s))'(x) \D s.$$
		
		Since~$u_0 \in \CC_0(\R,\R)$ is uniformly continuous, standard results imply that~$[0,\infty) \times \R \to \R$, $(t,x) \mapsto Q_t u_0(x)$ is uniformly continuous. Consequently, the first term in the definition of~$\tilde{u}$ has the desired continuity property. Regarding the remaining two terms, 
		it suffices by the triangle inequality to consider the two coordinate directions separately and establish the existence of uniform moduli of continuity. Let~$T \ge 0$. We start with continuity in space. Let~$t \in [0,T]$ as well as~$\delta \gne 0$ and~$x,z \in \R$ with~$\abs x-z\abs = \delta$. We have
		\begin{align*}
			&\left\abs \int_0^t \int_\R \e^{t-s} q_{t-s}(y-x) \D \beta_s(y) \D s - \int_0^t \int_\R \e^{t-s} q_{t-s}(y-z) \D \beta_s(y) \D s \right\abs\\
			\le{}& \e^T \int_0^t \sup_{\tilde{y} \in \R} \abs q_{t-s}(\tilde{y})- q_{t-s}(\tilde{y}+\delta)\abs \D s
			\le \e^T \int_0^t \min\{\delta \y q_{t-s}'\y_\infty, 2 \y q_{t-s}\y_\infty\} \D s
			\le \e^T \int_0^t \min\left\{\frac{\delta}{t-s}, \frac{2}{\sqrt{t-s}}\right\} \D s.
		\end{align*}	
		Noting
		\begin{align*}
			\int_\R \abs q'_{t-s}(\tilde{y}) - q'_{t-s}(\tilde{y}+\delta)\abs \D \tilde{y}
			\le{}& \min\{\delta \y q''_{t-s}\y_1, 2 \y q'_{t-s}\y_1\} \le \min\left\{ \frac{\delta}{t-s}, \frac{2}{\sqrt{t-s}}\right\}, \quad s \in (0,t],
		\end{align*}
		we also get
		\begin{align*}
			\left\abs \int_0^t \e^{t-s} (Q_{t-s} (bu_s))'(x) \D s - \int_0^t \e^{t-s} (Q_{t-s} (bu_s))'(z) \D s\right\abs
			={}& \left\abs \int_0^t \e^{t-s} \int_\R (q'_{t-s}(x-y) - q'_{t-s}(z-y)) b(y) u_s(y) \D y \D s\right\abs\\
			\le{}& \e^T \y b\y_\infty C_T \int_0^t \min\left\{\frac{\delta}{t-s}, \frac{2}{\sqrt{t-s}}\right\} \D s.
		\end{align*}
		It remains to observe that the dominated convergence theorem implies~$$\lim_{\delta \searrow 0} \sup_{t \in [0,T]} \int_0^t \min\left\{\frac{\delta}{t-s}, \frac{2}{\sqrt{t-s}}\right\} \D s = \lim_{\delta \searrow 0} \int_0^T \min\left\{\frac{\delta}{r}, \frac{2}{\sqrt{r}}\right\} \D r =0.$$
		We now turn toward continuity in space. Let~$x \in \R$ as well as~$\delta \in (0,T]$ and~$t \in [0,T-\delta]$. Noting~$\abs \partial_r q_r(y)\abs \le \frac{1}{\sqrt{r^3}}$ and thus, by the mean value theorem,
		\begin{align*}
			\abs \e^\delta q_{r+\delta}(y)-q_r(y)\abs \le (\e^\delta -1) q_{r+\delta}(y) + \abs q_{r+\delta}(y)-q_r(y)\abs 
			\le \min\left\{\frac{\e^\delta-1}{\sqrt{r}} + \frac{\delta}{\sqrt{r^3}} ,\frac{\e^\delta+1}{\sqrt{r}}\right\}, \quad r \gne 0,\, y \in \R,
		\end{align*}
		we get
		\begin{align*}
			&\left\abs \int_0^{t+\delta} \int_\R \e^{t+\delta-s} q_{t+\delta-s}(y-x) \D \beta_s(y) \D s - \int_0^t \int_\R \e^{t-s} q_{t-s}(y-x) \D \beta_s(y) \D s \right\abs\\
			\le{}& \e^T \int_t^{t+\delta} \frac{1}{\sqrt{t+\delta-s}} \D s + \e^T\int_0^t \sup_{y \in \R} \abs \e^\delta q_{t+\delta-s}(y)-q_{t-s}(y)\abs \D s\\
			\le{}& \e^T \int_0^\delta \frac{1}{\sqrt{r}} \D r + \e^T \int_0^T \min\left\{\frac{\e^\delta-1}{\sqrt{r}} + \frac{\delta}{\sqrt{r^3}} ,\frac{\e^\delta+1}{\sqrt{r}}\right\} \D r,
		\end{align*}
		which converges to~$0$ as~$\delta \searrow 0$ by the dominated convergence theorem. Using~$\y q_r'\y_1 \le \frac{1}{\sqrt{r}}$ and~$\y \partial_r q_r'\y_1 \le \frac{2}{\sqrt{r^3}}$, the mean value theorem again yields
		\begin{align*}
			\y\e^\delta q'_{r+\delta} - q'_r\y_1
			\le (\e^\delta-1) \y q'_{r+\delta}\y_1 + \y q'_{r+\delta} - q'_r\y_1
			\le \min\left\{\frac{\e^\delta-1}{\sqrt{r}}+\frac{2\delta}{\sqrt{r^3}},\frac{\e^\delta+1}{\sqrt{r}}\right\}, \quad r \gne 0.
		\end{align*}
		We finally obtain
		\begin{align*}
			&\left\abs \int_0^{t+\delta} \e^{t+\delta-s} (Q_{t+\delta-s} (bu_s))'(x) \D s - \int_0^t \e^{t-s} (Q_{t-s} (bu_s))'(x) \D s\right\abs\\
			={}& \left\abs \int_0^{t+\delta} \int_\R \e^{t+\delta-s} q'_{t+\delta-s}(x-y) b(y) u_s(y) \D y \D s - \int_0^t \int_\R \e^{t-s} q'_{t-s}(x-y) b(y) u_s(y) \D y \D s \right\abs\\
			\le{}& \e^T \y b\y_\infty C_T \left(\int_t^{t+\delta} \y q'_{t+\delta-s}\y_1 \D s + \int_0^t \y\e^\delta q'_{t+\delta-s} - q'_{t-s}\y_1 \D s\right)\\
			\le{}& \e^T \y b\y_\infty C_T \left(\int_0^\delta \frac{1}{\sqrt{r}} \D r + \int_0^T \min\left\{\frac{\e^\delta-1}{\sqrt{r}}+ \frac{2\delta}{\sqrt{r^3}}, \frac{\e^\delta+1}{\sqrt{r}}\right\} \D r \right)
		\end{align*}
		which again converges to~$0$ as~$\delta \searrow 0$.
	\end{proof}

	\section{Uniqueness of PDE Solution}
	\label{sec:unique}
	Our proof builds on the uniqueness argument of~\cite[Section 5]{AMW}. A key difficulty in extending that argument to discontinuous drifts
	is that its estimates rely on regularity of $b$. We overcome this
	difficulty by controlling the contribution of the drift using only
	its boundedness, and obtain the following result.
	\begin{prop}
		\label{prop:unique}
		In the space~$\UU \times \MM_{\nearrow 1}$, the PDE problem~\eqref{eqn:PDE} has at most one solution.
	\end{prop}
	
	Throughout this section, let~$t \gne 0$ be fixed, $\eta \in \CC_c^\infty(\R,[0,\infty))$ a symmetric mollifier supported on~$[-1,1]$, and~$\delta \gne 0$. Further, we use the following notation: Introducing~$\eta_\delta:\R \to [0,\infty)$, $x \mapsto \frac{1}{\delta}\eta(\frac{x}{\delta})$, we set~$\psi^\delta:=\eta_\delta * \psi$ for~$\psi: \R \to \R$. For~$\nu \in \MM_1$, we denote the smooth Lebesgue density of~$\eta_\delta * \nu$ by~$\nu^\delta$. Finally, we write~$w_s:=w(s,\cdot)$ for~$w \in \UU$ and~$s \in [0,t]$ to simplify the notation.
	
	We start with an auxiliary result stating that, for any~$u \in \UU$, certain mollifications in space converge locally uniformly in time.
	\begin{lem}
		\label{lem:moll}
		Let~$w \in \UU$. Then we have $$\lim_{\delta \searrow 0} \sup_{s \in [0,t]} \y w_s - (w_s^\delta)^\delta\y_\infty = 0 \quad \text{and} \quad \lim_{\delta \searrow 0} \sup_{s \in [0,t]} \y w_s - w_s^\delta\y_2 = 0.$$
	\end{lem}
	\begin{proof}
		Regarding the first claim, it is straightforward to check that~$w\abs_{[0,t]\times \R}$ is uniformly continuous: On the one hand, for every~$\eps \gne 0$, there exists an~$R \gne 0$ with~$\sup_{s \in [0,t]} \sup_{x \in [-R,R]^c} \abs w_s(x)\abs \lne \eps$. On the other hand, $w$ is uniformly continuous on the compact set~$[0,t] \times [-R,R]$. Noting~$\int_\R \eta_\delta(x) \D x = 1$ and~$\op{supp}(\eta_\delta) \subseteq [-\delta,\delta]$, we get
		\begin{align*}
			\abs w_s(x) - (w_s^\delta)^\delta(x)\abs 
			={}& \left\abs \int_\R \int_\R \eta_\delta(x-y) \eta_\delta(y-z) (w_s(x)-w_s(z)) \D z \D y\right\abs\\
			\le{}& \int_{x-\delta}^{x+\delta} \int_{y-\delta}^{y+\delta} \eta_\delta(x-y) \eta_\delta(y-z) \abs w_s(x)-w_s(z)\abs \D z \D y\\
			\le{}& \sup_{z \in [x-2\delta,x+2\delta]} \abs w_s(x)-w_s(z)\abs, \qquad x \in \R,\, s \in [0,t].
		\end{align*}
		As~$\delta \searrow 0$, this converges to~$0$ uniformly in~$x \in \R$ and~$s \in [0,t]$.
		
		Let us now turn toward the second claim. Let~$\eps \gne 0$. As a consequence of~$w \in \UU$, there exists an~$R \gne 1$ with~$\sup_{s \in [0,t]} \int_{[-(R-1),R-1]^c} \abs w_s(x)\abs^2 \D x \lne \frac{\eps}{8}$. We obtain
		\begin{align*}
			&\int_{[-R,R]^c} \int_{-\delta}^\delta \abs w_s(x-y)-w_s(x)\abs^2\eta_\delta(y) \D y \D x\\
			\le{}& 2 \int_{[-R,R]^c} \int_{-\delta}^\delta \abs w_s(x-y)\abs^2 \eta_\delta(y) \D y \D x + 2 \int_{[-R,R]^c} \int_{-\delta}^\delta \abs w_s(x)\abs^2\eta_\delta(y) \D y \D x\\
			\le{}& 4 \int_{[-(R-1),R-1]^c} \abs w_s(x)\abs^2 \D x \int_{-\delta}^\delta \eta_\delta(y) \D y
			\lne \frac{\eps}{2}, \qquad s \in [0,t],\, \delta \in (0,1].
		\end{align*}
		Since~$w\abs_{[0,t] \times [-(R+1),R+1]}$ is uniformly continuous, there exists a~$\delta_0 \in (0,1]$ with~$\abs w_s(x)-w_s(y)\abs \lne \sqrt{\frac{\eps}{4R}}$ for all~$s \in [0,t]$ and all~$x,y \in [-(R+1),R+1]$ with~$\abs x-y\abs \lne \delta_0$. We deduce
		\begin{align*}
			\int_{-R}^R \int_{-\delta}^\delta \abs w_s(x-y)-w_s(x)\abs^2\eta_\delta(y) \D y \D x
			\lne{}& \int_{-R}^R \frac{\eps}{4R} \D x \int_{-\delta}^\delta \eta_\delta(y) \D y = \frac{\eps}{2}, \quad s \in [0,t],\, \delta \in (0,\delta_0].
		\end{align*}
		Combining Jensen's inequality, applied to the probability measure with density~$\eta_\delta$, with the two bounds we just established, we finally get
		\begin{align*}
			\y w_s - w_s^\delta\y_2^2
			={}& \int_\R \left\abs \int_{-\delta}^\delta (w_s(x-y)-w_s(x))\eta_\delta(y) \D y\right\abs^2 \D x\\
			\le{}& \int_\R \int_{-\delta}^\delta \abs w_s(x-y)-w_s(x)\abs^2\eta_\delta(y) \D y \D x
			\lne \eps, \qquad s \in [0,t],\,\delta \in (0,\delta_0],
		\end{align*}
		proving~$\lim_{\delta \searrow 0} \sup_{s \in [0,t]} \y w_s - w_s^\delta\y_2 = 0$.
	\end{proof}

	\begin{proof}[Proof of Proposition~\ref{prop:unique}]
		Let~$(u,\beta),(\tilde{u},\tilde{\beta}) \in \UU \times \MM_{\nearrow 1}$ be two solutions to~\eqref{eqn:PDE}. Let us start by briefly arguing that it suffices to show~$u=\tilde{u}$. Indeed, by definition of a solution, $u=\tilde{u}$ implies $$\int_{[0,t]\times \R} \varphi(x) \D \beta(s,x) = \int_{[0,t]\times \R} \varphi(x) \D \tilde{\beta}(s,x), \quad t \ge 0,\, \varphi \in \CC_c^\infty(\R,\R).$$ Recalling the definition of~$\MM_{\nearrow 1}$ and invoking standard arguments from measure theory, this, in turn, yields~$\beta = \tilde{\beta}$.
		
		Now let~$v:=u-\tilde{u} \in \UU$ and~$\varphi \in \CC_c^\infty(\R,\R)$. The definition of a solution, applied to~$\varphi^\delta \in \CC_c^\infty(\R,\R)$, yields
		\begin{align*}
			\int_\R \varphi^\delta (u_t - u_0) \D \lambda
			={}& \int_0^t \int_\R \left(\frac{1}{2} (\varphi^\delta)'' + (\varphi^\delta)' b + \varphi^\delta \right) u_s \D \lambda \D s
			- \int_0^t \int_\R \varphi^\delta \D \beta_s \D s
		\end{align*}
		Noting~$v_0=u_0-u_0=0$, subtracting the analogous equation for~$(\tilde{u},\tilde{\beta})$ yields
		\begin{align*}
			\int_\R \varphi^\delta v_t \D \lambda
			={}& \int_0^t \int_\R \left(\frac{1}{2} (\varphi^\delta)'' + (\varphi^\delta)' b + \varphi^\delta \right) v_s \D \lambda \D s
			- \int_0^t \int_\R \varphi^\delta \D \beta_s \D s + \int_0^t \int_\R \varphi^\delta \D \tilde{\beta}_s \D s.
		\end{align*}
		By antisymmetry of~$\eta_\delta'$ and Fubini's theorem, we have
		\begin{align}
			\label{eqn:antisym}
			\int_\R (\varphi^\delta)' b  v_s \D \lambda
			={}& \int_\R \int_\R \eta_\delta'(x-y) \varphi(y) b(x) v_s(x) \D y \D x \notag\\
			={}& -\int_\R \int_\R \eta_\delta'(y-x) \varphi(y) b(x) v_s(x) \D x \D y
			= -\int_\R \varphi \cdot ((bv_s)^\delta)' \D \lambda.
		\end{align}
		Together with similar arguments for the remaining terms, the above equation can be rewritten as
		\begin{align*}
			\int_\R \varphi v_t^\delta \D \lambda
			={}& \int_0^t \int_\R \varphi \left(\frac{1}{2} (v_s^\delta)'' - ((bv_s)^\delta)' + v_s^\delta \right) \D \lambda \D s
			+ \int_0^t \int_\R \varphi (\tilde{\beta}_s^\delta-\beta_s^\delta) \D \lambda \D s.
		\end{align*}
		Since this holds for every~$\varphi \in \CC_c^\infty(\R,\R)$, the fundamental lemma of the calculus of variations yields
		\begin{align*}
			v_t^\delta
			={}& \int_0^t \left(\frac{1}{2} (v_s^\delta)'' - ((bv_s)^\delta)' + v_s^\delta \right) \D s
			+ \int_0^t \tilde{\beta}_s^\delta-\beta_s^\delta \D s.
		\end{align*}
		Before continuing, we record a few observations: We have~$v \in \UU$, which guarantees~$\sup_{s \in [0,t]} \y v_s\y_2 \lne \infty$. Applying Young's convolution inequality and recalling that~$b$ is bounded, we deduce
		\begin{align*}
			\sup_{s \in [0,t]} \y v_s^\delta\y_2 \lne \infty,
			\quad \sup_{s \in [0,t]} \y (v_s^\delta)'' \y_2 \lne \infty, \quad \sup_{s \in [0,t]} \y ((bv_s)^\delta)' \y_2 \lne \infty.
		\end{align*}
		Moreover, $v \in \UU$ guarantees~$\sup_{s \in [0,t]} \y v_s^\delta\y_\infty \le \sup_{s \in [0,t]} \y v_s\y_\infty \lne \infty$. Finally, $\beta_s^\delta$ and~$\tilde{\beta}_s^\delta$ are probability densities for all~$s \in [0,t]$ since~$\beta_s$ and~$\tilde{\beta}_s$ are probability measures.
		These observations, together with Hölder's inequality, allow us to apply Fubini's theorem in the subsequent calculation. Noting~$v_0^\delta = 0$, the above representation of~$v_t^\delta$ yields
		\begin{align}
			\label{eqn:Gw0}
			\frac{1}{2} \y v_t^\delta \y_2^2 = \int_\R \int_0^t v_s^\delta \dot{v}^\delta_s \D s \D \lambda
			= \int_0^t \int_\R v_s^\delta \left(\frac{1}{2} (v_s^\delta)'' -  ((bv_s)^\delta)' + v_s^\delta \right) \D \lambda \D s + \int_0^t \int_\R v_s^\delta (\tilde{\beta}_s^\delta-\beta_s^\delta) \D \lambda \D s.
		\end{align}
		Our aim now is to prepare an application of Gronwall's lemma. To this end, we consider the terms on the right-hand side of~\eqref{eqn:Gw0} one by one.
		Since~$v \in \UU$ implies~$v_s \in \CC_0(\R,\R)$ and hence~$v_s^\delta \in \CC_0(\R,\R)$ for all~$s \in [0,t]$, integration by parts yields
		\begin{align}
			\label{eqn:Gw1}
			\int_0^t \int_\R v_s^\delta (v_s^\delta)'' \D \lambda \D s
			= - \int_0^t \int_\R (v_s^\delta)' (v_s^\delta)' \D \lambda \D s
			= - \int_0^t \y (v_s^\delta)'\y_2^2 \D s.
		\end{align}
		Combining a similar integration by parts with the Cauchy-Schwarz inequality, we obtain
		\begin{align*}
			- \int_0^t \int_\R v_s^\delta ((bv_s)^\delta)' \D \lambda \D s
			={}& \int_0^t \int_\R (v_s^\delta)' (bv_s)^\delta \D \lambda \D s\\
			\le{}& \int_0^t \y (v_s^\delta)'\y_2 \y (bv_s)^\delta\y_2 \D s
			\le \frac{1}{2} \int_0^t \y (v_s^\delta)'\y_2^2 \D s + \frac{1}{2} \int_0^t\y (bv_s)^\delta\y_2^2 \D s.
		\end{align*}
		Further, we get
		\begin{align*}
			\y (bv_s)^\delta\y_2^2
			={}& \int_\R \left\abs \int_\R \eta_\delta(x-y) b(y) v_s(y) \D y \right\abs^2 \D x
			\le \int_\R \left(\int_\R \eta_\delta(x-y) \abs b(y)\abs  \abs v_s(y)\abs  \D y\right)^2 \D x\\
			\le{}& \y b\y_\infty^2 \int_\R \left(\int_\R \eta_\delta(x-y) \abs v_s(y)\abs  \D y\right)^2 \D x
			= \y b\y_\infty^2 \y \abs v_s\abs^\delta\y_2^2, \qquad s \in [0,t].
		\end{align*}
		Noting~$\y v_s\y_2 = \y \abs v_s\abs\y_2$, we further obtain
		\begin{align*}
			\y \abs v_s\abs^\delta\y_2^2
			\le{}& (\y \abs v_s\abs^\delta - \abs v_s\abs \y_2 + \y v_s - v_s^\delta\y_2 + \y v_s^\delta\y_2)^2
			\le 2\big(\big\y \abs v_s\abs^\delta - \abs v_s\abs \big\y_2 + \big\y v_s - v_s^\delta\big\y_2\big)^2 + 2\y v_s^\delta\y_2^2, \quad s \in [0,t].
		\end{align*}
		Combining these three inequalities, we deduce
		\begin{align}
			\label{eqn:Gw2}
			&- \int_0^t \int_\R v_s^\delta ((bv_s)^\delta)' \D \lambda \D s\notag\\
			\le{}& \frac{1}{2} \int_0^t \y (v_s^\delta)'\y_2^2 \D s + \y b\y_\infty^2 t \sup_{s \in [0,t]} \big(\big\y \abs v_s\abs^\delta - \abs v_s\abs \big\y_2 + \big\y v_s - v_s^\delta\big\y_2\big)^2 + \y b\y_\infty^2 \int_0^t\y v_s^\delta\y_2^2 \D s.
		\end{align}
		Let us finally consider the last term in~\eqref{eqn:Gw0}. Recalling~$\beta \in \BB^u$ and~$u \ge 0$, we have~$\int_\R u_s \D \beta_s = \y u_s\y_\infty \ge \int_\R u_s \D \tilde{\beta}_s$ for a.e.~$s \in [0,t]$. Together with the analogous fact for~$\tilde{u}$, we get~$\int_\R v_s \D (\tilde{\beta}_s-\beta_s)\le 0$ for a.e.~$s \in [0,t]$. Using an argument similar to~\eqref{eqn:antisym}, based on the symmetry of~$\eta_\delta$, in the first step, we thus obtain 
		\begin{align*}
			\int_0^t \int_\R v_s^\delta (\tilde{\beta}_s^\delta-\beta_s^\delta) \D \lambda \D s
			={}& \int_0^t \int_\R (v_s^\delta)^\delta \D (\tilde{\beta}_s- \beta_s) \D s\\
			\le{}& \int_0^t \int_\R v_s \D (\tilde{\beta}_s-\beta_s) \D s + 2 \int_0^t \y v_s - (v_s^\delta)^\delta\y_\infty \D s
			\le 2t \sup_{s \in [0,t]} \y v_s - (v_s^\delta)^\delta\y_\infty.
		\end{align*}
		Inserting this as well as~\eqref{eqn:Gw1} and~\eqref{eqn:Gw2} into~\eqref{eqn:Gw0}, we get
		\begin{align*}
			\frac{1}{2} \y v_t^\delta \y_2^2 \le (\y b\y_\infty^2+1) \int_0^t \y v_s^\delta\y_2^2 \D s + \y b\y_\infty^2 t \sup_{s \in [0,t]} (\y \abs v_s\abs^\delta - \abs v_s\abs \y_2 + \y v_s - v_s^\delta\y_2)^2 + 2t \sup_{s \in [0,t]} \y v_s - (v_s^\delta)^\delta\y_\infty.
		\end{align*}
		Gronwall's lemma yields
		\begin{align*}
			\y v_t^\delta\y_2^2 \le \left(2\y b\y_\infty^2 t \sup_{s \in [0,t]} (\y \abs v_s\abs^\delta - \abs v_s\abs \y_2 + \y v_s - v_s^\delta\y_2)^2 + 4t \sup_{s \in [0,t]} \y v_s - (v_s^\delta)^\delta\y_\infty \right) \e^{2(\y b\y_\infty^2+1)t}.
		\end{align*}
		Noting~$v, \abs v \abs \in \UU$, Lemma~\ref{lem:moll} shows that this converges to~$0$ as~$\delta \searrow 0$, proving~$\y v_t\y_2 = \lim_{\delta \searrow 0} \y v_t^\delta\y_2 = 0$. Since~$v_t$ is continuous, we end up with~$u(t,\cdot\,) = \tilde{u}(t,\cdot\,)$.
	\end{proof}

	\section{Complementarity Condition and Proof of Theorem~\ref{thm:1}}
	\label{sec:final}
	
	One of our main tasks in this final section can roughly be summarized as follows: If~$(\mu,\beta)$ is a subsequential distributional limit of~$(\mu^{n,\eps},\beta^{n,\eps})$ as~$n \to \infty$ followed by~$\eps \searrow 0$ and if~$u$ is a continuous density of~$\mu$, then the complementarity condition~$\beta \in \BB^u$ is satisfied a.s. The argument for the limit~$\eps \searrow 0$ will be based on uniform convergence of the densities on every finite time horizon.
	
	\begin{lem}
		\label{lem:AA}
		Let~$(\mu^m)_{m \in \N} \subseteq \CC([0,\infty),\MM_1)$ be a sequence which converges to~$\mu \in \CC([0,\infty),\MM_1)$. We assume that~$\mu^m$ has a density~$u^m$ for every~$m \in \N$ such that~$(u^m\abs_{[0,T] \times \R})_{m \in \N}$ is uniformly equicontinuous for every~$T \ge 0$. For each~$T \ge 0$, we additionally assume\footnote{These additional assumptions are actually implied by the properties stated before. However, since we already know that the additional assumptions are satisfied in our setting, we simply include them to shorten the proof.} $\sup_{m \in \N} \y u^m\abs_{[0,T] \times \R}\y_\infty \lne \infty$ and
		\begin{align}
			\label{eqn:decay1}
			\lim_{R \to \infty} \sup_{m \in \N} \sup_{t \in [0,T]} \sup_{x \in [-R,R]^c} \abs u^m(t,x)\abs = 0.
		\end{align}
		Then~$\mu$ has a density~$u \in \CC([0,\infty) \times \R,\R)$ satisfying
		\begin{align}
			\label{eqn:decay2}
			\lim_{m \to \infty} \y (u^m-u)\abs_{[0,T] \times \R}\y_\infty = 0 \quad \text{and} \quad \lim_{R \to \infty} \sup_{t \in [0,T]} \sup_{x \in [-R,R]^c} \abs u(t,x)\abs = 0, \quad T \ge 0.
		\end{align}
		Moreover, if we also have
		\begin{align}
			\label{eqn:decay3}
			\lim_{R \to \infty} \sup_{m \in \N} \sup_{t \in [0,T]} \int_{[-R,R]^c} \abs u^m(t,x)\abs^2 \D x = 0,
		\end{align}
		then~$u \in \UU$ holds.
	\end{lem}
	\begin{proof}
		The Arzelà-Ascoli theorem (combined with exhausting~$[0,\infty) \times \R$ by a sequence of compact sets and a diagonalization argument) implies that every subsequence of~$(u^m)_{m \in \N}$ has a sub-subsequence converging locally uniformly to some~$u \in \CC([0,\infty) \times \R,\R)$. Now let~$t \ge 0$. On the one hand, convergence of~$(\mu^m)_{m \in \N}$ to~$\mu$ yields~$$\int_\R \varphi(x) u^m(t,x) \D x = \int_\R \varphi(x) \D \mu_t^m(x) \xrightarrow{m \to \infty} \int_\R \varphi(x) \D \mu_t(x), \quad \varphi \in \CC_c(\R,\R).$$ On the other hand, the left-hand side converges to~$\int_\R \varphi(x) u(t,x) \D x$ for every~$\varphi \in \CC_c(\R,\R)$ along the sub-subsequence from above. Consequently, $u$ is a density of~$\mu$. Since~$\mu$ can have at most one continuous density, the whole sequence~$(u^m)_{m \in \N}$ must converge locally uniformly to~$u$ by a standard sub-subsequence argument.
		Together with~\eqref{eqn:decay1}, we obtain~\eqref{eqn:decay2}. The final claim is an immediate consequence of the convergence just established, Fatou's lemma and the definition of~$\UU$.
	\end{proof}
	
	We proceed with the announced auxiliary result concerning the complementarity condition.
	\begin{lem}
		\label{lem:comp}
		For every~$m \in \N$, let~$(\mu^m,\beta^m)$ be any subsequential distributional limit of~$(\mu^{n,1/m},\beta^{n,1/m})_{n \in \N}$. Further, let~$(\mu,\beta)$ be any subsequential distributional limit of~$(\mu^m,\beta^m)_{m \in \N}$ and let~$u$ a.s.~be a (random) continuous density of~$\mu$. Then~$\beta \in \BB^u$ holds a.s.
	\end{lem}
	\begin{proof}
		Throughout the proof, let~$T \gne 0$ be fixed. For all~$n,m \in \N$, we introduce $$v^{n,m}: [0,T] \times \R \to [0,1], \quad (t,x) \mapsto
		\begin{cases}
			\mu_{t-}^{n,1/m}(C_k^{1/m}), \quad& x \in (\frac{k}{m},\frac{k+1}{m}) \text{ for some } k \in \Z,\\
			\max\{\mu_{t-}^{n,1/m}(C_{k-1}^{1/m}),\mu_{t-}^{n,1/m}(C_k^{1/m})\}, \quad& x=\frac{k}{m} \hspace{28pt} \text{ for some } k \in \Z,
		\end{cases}$$
		with the convention~$\mu_{0-}^{n,1/m}:=\mu_0^{n,1/m}$. Noting~$$\Z \setminus A^{1/m}(X^{n,1/m}(t-)) = \left\{k \in \Z: \mu_{t-}^{n,1/m}(C_k^{1/m}) \lne \max_{k' \in \Z} \mu_{t-}^{n,1/m}(C_{k'}^{1/m})\right\}, \quad t \gne 0,\, n,m \in \N,$$
		we, by construction, have
		\begin{align*}
			&\beta^{n,1/m}(\{(t,x) \in [0,T] \times \R: v^{n,m}(t,x) \lne \y v^{n,m}(t,\cdot)\y_\infty \})\\
			\le{}&\beta^{n,1/m}(\{0\}\times \R) + \beta^{n,1/m}\left(\left\{(t,x) \in (0,T] \times \R: x \in \bigcup_{k \in \Z \setminus A^{1/m}(X^{n,1/m}(t-))} C_k^{1/m}\right\}\right)=0, \quad n,m \in \N,
		\end{align*}
		and thus
		\begin{align}
			\label{eqn:comp0}
			\int_{[0,T] \times \R} (\y v^{n,m}(t,\cdot)\y_\infty - v^{n,m}(t,x)) \D \beta^{n,1/m}(t,x) = 0, \quad n,m \in \N.
		\end{align}	
		
		Now fix~$m \in \N$. Our intermediate goal is to prove a similar concentration of~$\beta^m$. Let~$(n_l^m)_{l \in \N} \subseteq \N$ be a sequence of indices such that~$(\mu^m,\beta^m)$ is the distributional limit of~$(\mu^{n_l^m,1/m},\beta^{n_l^m,1/m})_{l \in \N}$. For the rest of this paragraph, we assume sure convergence by invoking Skorokhod's representation theorem and discarding the exceptional set.	All statements are to be understood pointwise on the new sample space. As~$l \to \infty$, we claim that~$v^{n_l^m,m}$ converges uniformly to
		$$v^m: [0,T] \times \R \to [0,1], \quad (t,x) \mapsto
		\begin{cases}
			\mu_t^m(C_k^{1/m}), \quad& x \in (\frac{k}{m},\frac{k+1}{m}) \text{ for some } k \in \Z,\\
			\max\{\mu_t^m(C_{k-1}^{1/m}),\mu_t^m(C_k^{1/m})\}, \quad& x=\frac{k}{m} \hspace{28pt} \text{ for some } k \in \Z.
		\end{cases}$$
		To prove this, let~$d_{\op{LP}}$ be the Lévy-Prohorov metric on~$\MM_1$. Since Proposition~\ref{prop:tight} guarantees~$\mu^m \in \CC([0,\infty),\MM_1)$, a classic property of Skorokhod spaces yields~$\lim_{l \to \infty} \sup_{t \in [0,T]} d_{\op{LP}}(\mu^{n_l^m,1/m}_t,\mu^m_t) = 0$. Noting that~$\mu^m\abs_{[0,T]}$ is uniformly continuous, this implies that~$(\eta_l)_{l \in \N} := (\sup_{t \in [0,T]} d_{\op{LP}}(\mu^{n_l^m,1/m}_{t-}, \mu^m_t))_{l \in \N}$ is a null sequence, essentially by means of the triangle inequality. Since~$\mu^m$ has a density which is bounded on~$[0,T] \times \R$ by Lemma~\ref{lem:density}, we have~$\lim_{\eta \searrow 0} \sup_{t \in [0,T]} \sup_{x \in \R} \mu^m_t((x,x+3\eta))=0$. Now let~$l \in \N$. By definition of~$d_{\op{LP}}$ and~$C_k^m$, we get
		\begin{align*}
			\mu^m_t(C_k^m)
			\le{}& \mu^{n_l^m,1/m}_{t-}(C_k^m) + \mu^{n_l^m,1/m}_{t-}((\tfrac{k}{m}-\eta_l,\tfrac{k}{m})) + \mu^{n_l^m,1/m}_{t-}([\tfrac{k+1}{m},\tfrac{k+1}{m}+\eta_l)) + \eta_l\\
			\le{}& \mu^{n_l^m,1/m}_{t-}(C_k^m) + \mu^m_t((\tfrac{k}{m}-2\eta_l,\tfrac{k}{m}+\eta_l)) + \mu^m_t((\tfrac{k+1}{m}-\eta_l,\tfrac{k+1}{m}+2\eta_l)) + 3\eta_l, \qquad t \in [0,T],\, k \in \Z,
		\end{align*}
		and
		\begin{align*}
			\mu^{n_l^m,1/m}_{t-}(C_k^m)
			\le{}& \mu^m_t(C_k^m) + \mu^m_t((\tfrac{k}{m}-\eta_l,\tfrac{k}{m})) + \mu^m_t([\tfrac{k+1}{m},\tfrac{k+1}{m}+\eta_l)) + \eta_l, \quad t \in [0,T],\, k \in \Z.
		\end{align*}
		Recalling~$\lim_{l \to \infty} \eta_l = 0$, we thus obtain
		\begin{align*}
			\limsup_{l \to \infty} \sup_{t \in [0,T]} \sup_{k \in \Z} \big\abs \mu^{n_l^m,1/m}_{t-}(C_k^m) - \mu^m_t(C_k^m)\big\abs
			\le{}& \lim_{l \to \infty} \left(2 \sup_{t \in [0,T]} \sup_{x \in \R} \mu_t^m((x,x+3\eta_l)) + 3\eta_l\right)
			= 0.
		\end{align*}
		This implies that~$v^{n_l^m,m}$ converges uniformly to~$v^m$ as~$l \to \infty$, as claimed. Together with~\eqref{eqn:comp0} and~$\beta^{n_l^m,1/m} \in \MM_{\nearrow 1}$ for all~$l \in \N$, we deduce
		\begin{align*}
			&\lim_{l \to \infty} \int_{[0,T] \times \R} (\y v^m(t,\cdot)\y_\infty - v^m(t,x)) \D \beta^{n_l^m,1/m}(t,x)\\
			={}&\lim_{l \to \infty} \int_{[0,T] \times \R} \big((\y v^m(t,\cdot)\y_\infty - v^m(t,x)) - (\y v^{n_l^m,m}(t,\cdot)\y_\infty - v^{n_l^m,m}(t,x))\big) \D \beta^{n_l^m,1/m}(t,x) = 0.
		\end{align*}
		Finally, since the integrand is lower semicontinuous and non-negative, the Portmanteau theorem yields
		\begin{align}
			\label{eqn:comp1}
			\int_{[0,T] \times \R} (\y v^m(t,\cdot)\y_\infty - v^m(t,x)) \D \beta^m(t,x) = 0 \quad \text{a.s.}
		\end{align}
		
		Let us now turn toward the limit as~$m \to \infty$. As above, we invoke Skorokhod's representation theorem to assume that~$(\mu^m,\beta^m)_{m \in \N}$ surely converges to~$(\mu,\beta)$ along a subsequence. Lemmas~\ref{lem:density}, \ref{lem:PDE} and~\ref{lem:cont} imply that~$\mu^m$ has a (random) continuous density~$u^m$ for each~$m \in \N$ such that~$(u^m\abs_{[0,T] \times \R})_{m \in \N}$ is uniformly equicontinuous and uniformly bounded with~\eqref{eqn:decay1}. Thus Lemma~\ref{lem:AA} shows that~$\mu$ a.s.~has a (random) density~$\tilde{u} \in \CC([0,\infty) \times \R, \R)$ such that~$(u^m\abs_{[0,T] \times \R})_{m \in \N}$ converges uniformly to~$\tilde{u}\abs_{[0,T] \times \R}$ along the subsequence along which~$(\mu^m,\beta^m)_{m \in \N}$ converges to~$(\mu,\beta)$. Since~$\mu$ can have at most one continuous density, we must a.s.~have~$\tilde{u}=u$. Now let~$w: (0,\infty) \to [0,\infty)$ be a joint modulus of continuity for~$(u^m\abs_{[0,T] \times \R})_{m \in \N}$. By construction of~$v^m$ and by the mean value theorem for integrals, we a.s.~obtain~$\y u^m\abs_{[0,T] \times \R} - m v^m\y_\infty \le w(\frac{1}{m})$ for all~$m \in \N$, showing that~$(u^m\abs_{[0,T] \times \R} - m v^m)_{m \in \N}$ a.s.~converges to~$0$. Together with~\eqref{eqn:comp1} and~$\beta^m \in \MM_{\nearrow 1}$ for all~$m \in \N$, we deduce that
		\begin{align*}
			&\int_{[0,T] \times \R} (\y u(t,\cdot)\y_\infty - u(t,x)) \D \beta^m(t,x)\\
			\le{}& \int_{[0,T] \times \R} \big((\y u(t,\cdot)\y_\infty - u(t,x)) - (\y u^m(t,\cdot)\y_\infty - u^m(t,x))\big) \D \beta^m(t,x)\\
			&\quad+ \int_{[0,T] \times \R} \left((\y u^m(t,\cdot)\y_\infty - u^m(t,x)) - m(\y v^m(t,\cdot)\y_\infty - v^m(t,x))\right) \D \beta^m(t,x)
		\end{align*}
		a.s. converges to~$0$ as~$m \to \infty$ along the subsequence from above. Once again, the Portmanteau theorem yields 
		\begin{align*}
			\int_{[0,T] \times \R} (\y u(t,\cdot)\y_\infty - u(t,x)) \D \beta(t,x) = 0 \quad \text{a.s.}
		\end{align*}
		Since~$T \gne 0$ was arbitrary, we must have $$\beta(\{(t,x) \in [0,\infty) \times \R: u(t,x) \lne \y u(t, \cdot\,)\y_\infty\})=0 \quad \text{a.s.},$$ as claimed.
	\end{proof}
	
	We continue with a non-standard version of a classic diagonalization result for limits of double sequences,
	which will be the top-level ingredient of our proof of Theorem~\ref{thm:1}.
	
	\begin{lem}
		\label{lem:diag}
		Let~$(\EE,d)$ be a compact metric space and let~$(z_{n,m})_{n,m \in \N} \subseteq \EE$. Assume that there is some~$z \in \EE$ with the following property: If~$z_m \in \EE$ is any accumulation point of~$(z_{n,m})_{n \in \N}$ for every~$m \in \N$, the sequence~$(z_m)_{m \in \N}$ converges to~$z$. Then there exists a sequence~$(m_n)_{n \in \N} \subseteq \N$ with~$\lim_{n \to \infty} m_n = \infty$ such that~$(z_{n,m_n})_{n \in \N}$ converges to~$z$.
	\end{lem}
	\begin{proof}
		Let~$l \in \N$. By assumption, there must be an~$m'_l \in \N$ such that, for every~$m \in \N$ with~$m \ge m'_l$, every accumulation point~$z_m \in \EE$ of~$(z_{n,m})_{n \in \N}$ satisfies~$d(z_m,z) \lne \frac{1}{2l}$. Due to compactness of~$\EE$, this implies that, for every~$m \in \N$ with~$m \ge m'_l$, there must be an~$n_{l,m} \in \N$ with~$d(z_{n,m},z) \lne \frac{1}{l}$ for all~$n \in \N$ with~$n \ge n_{l,m}$. Using the conventions~$\max\emptyset:=0$ and~$[0]:=\emptyset$, we set $$n_m:= m + \max\{n_{l,j}: l,j \in [m],\, j \ge m_l'\}, \quad m \in \N,$$ and define $$m_n:= m, \quad n \in [n_m] \setminus [n_{m-1}],\, m \in \N.$$ Since~$(n_m)_{m \in \N}$ is strictly increasing, $(m_n)_{n \in \N} \subseteq \N$ is a well-defined non-decreasing sequence with~$\lim_{n \to \infty} m_n = \infty$. Now let~$l \in \N$ and~$\bar{m} := \max\{l,m_l'\}$. Further, let~$n \in \N$ with~$n \ge n_{\bar{m}}$. Then we have~$m_n \ge \bar{m} \ge m_l'$ and thus~$n \ge n_{\bar{m}} \ge n_{l,m_n}$, proving~$d(z_{n,m_n},z) \lne \frac{1}{l}$.
	\end{proof}
	
	Finally, we are ready to put everything together.
	\begin{proof}[Proof of Theorem~\ref{thm:1}]
		For a Polish space~$E$, let~$\MM_1(E)$ be the Polish space of all Borel probability measures on~$E$. Further, for a random variable~$Z$, let~$\LL(Z)$ be its law. As a consequence of Proposition~\ref{prop:tight} and Prohorov's theorem, the metrizable space $$\EE:=\overline{\{\LL(\mu^{n,\eps},\beta^{n,\eps}): n \in \N,\,\eps \gne 0\}} \subseteq \MM_1(\DD([0,\infty),\MM_1) \times \MM_{\nearrow 1})$$ is compact. Now let~$m \in \N$ and let~$(\mu^m,\beta^m)$ be any subsequential distributional limit of~$(\mu^{n,1/m},\beta^{n,1/m})_{n \in \N}$. (Note that this is equivalent to~$(\mu^m,\beta^m) \sim \Q_m$ for some accumulation point~$\Q_m \in \EE$ of~$(\LL(\mu^{n,1/m},\beta^{n,1/m}))_{n \in \N}$.)
		Moreover, let~$(\tilde{\mu},\tilde{\beta})$ be a subsequential distributional limit of~$(\mu^m,\beta^m)_{m \in \N}$. Invoking Skorokhod's representation theorem, we assume a.s.~convergence. Lemmas~\ref{lem:density}, \ref{lem:PDE} and~\ref{lem:cont} imply that~$\mu^m$ has a (random) continuous density~$u^m$ with
		\begin{align}
			\label{eqn:PDEm}
			&\int_\R \varphi(x) u^m(t,x) \D x\notag\\
			={}& \int_\R \varphi(x) u_0(x) \D x
			- \int_{[0,t]\times \R} \varphi(x) \D \beta^m(s,x)\notag\\
			&\quad+ \int_0^t \int_\R \left(\frac{1}{2} \varphi''(x) + \varphi'(x) b(x) + \varphi(x) \right) u^m(s,x) \D x \D s, \qquad \varphi \in \CC_c^\infty(\R,\R),\, t \ge 0, \quad \text{a.s.},
		\end{align}
		and that~$(u^m\abs_{[0,T] \times \R})_{m \in \N}$ is a.s.~uniformly equicontinuous and uniformly bounded with~\eqref{eqn:decay1} and~\eqref{eqn:decay3} for all~$T \ge 0$. Since Proposition~\ref{prop:tight} a.s.~yields~$\mu^m \in \CC([0,\infty),\MM_1)$, we a.s.~have~$\tilde{\mu} \in \CC([0,\infty),\MM_1)$ as well. Lemma~\ref{lem:AA} a.s.~shows that~$\tilde{\mu}$ has a (random) density~$\tilde{u} \in \UU$ with~\eqref{eqn:decay2}. Taking the subsequential limit~$m \to \infty$ in~\eqref{eqn:PDEm}, using the dominated convergence theorem, we obtain 
		\begin{align*}
			\int_\R \varphi(x) \tilde{u}(t,x) \D x
			={}& \int_\R \varphi(x) u_0(x) \D x
			- \int_{[0,t]\times \R} \varphi(x) \D \tilde{\beta}(s,x)\\
			&\quad+ \int_0^t \int_\R \left(\frac{1}{2} \varphi''(x) + \varphi'(x) b(x) + \varphi(x) \right) \tilde{u}(s,x) \D x \D s, \qquad \varphi \in \CC_c^\infty(\R,\R),\, t \ge 0, \quad \text{a.s.}\notag
		\end{align*}
		Together with Lemma~\ref{lem:comp}, this shows that~$(\tilde{u},\tilde{\beta}) \in\UU\times\MM_{\nearrow1}$ solves~\eqref{eqn:PDE} a.s. In particular, a solution exists, and Proposition~\ref{prop:unique} ensures that it is unique.
		Denoting this solution by $(u,\beta)$, uniqueness gives~$(\tilde{u},\tilde{\beta})=(u,\beta)$ a.s. Letting~$\mu \in \CC([0,\infty),\MM_1)$ with density~$u$, this also implies~$\tilde{\mu}=\mu$ a.s. Recalling that convergence in distribution to a deterministic limit implies convergence in probability, Lemma~\ref{lem:diag} finally yields the claim.
	\end{proof}
	
	\section*{Acknowledgment}
	RA and DS were supported in part by ISF 3240/25. LM and DS were supported in part by ISF 2917/26. DS was supported in part by a Technion fellowship.
	
	\bibliographystyle{alpha}

\end{document}